\RequirePackage{etoolbox}
\csdef{input@path}{{style/}{graphics/}}
\documentclass[aap,MSNbibl,seceqn,dvips]{arximspdf}
\usepackage{mathrsfs}

\doi{10.1214/14-AAP1013} %kopijuoti is PTS
\volume{25}
\issue{2}
\pubyear{2015}
\firstpage{898}
\lastpage{929}

\makeatletter
\newtheorem{Thm}{Theorem}[section]
\newtheorem{Lem}[Thm]{Lemma}
\newproclaim{definition}[Thm]{Definition}

\newproclaim{remark}[Thm]{Remark}
\newproclaim{Example}{Example}
\newproclaim{Note}{Note}

\newcommand{\R}{\mathbb R}
\newcommand{\N}{\mathbb N}
\newcommand{\Z}{\mathbb Z}

\newcommand{\La}{\Lambda}
\newcommand{\De}{\Delta}
\newcommand{\al}{\alpha}
\newcommand{\be}{\beta}
\newcommand{\ga}{\gamma}
\newcommand{\si}{\sigma}
\newcommand{\ep}{\varepsilon}

\newcommand{\Om}{\Omega}

\newcommand{\F}{\mathcal{F}}
\newcommand{\G}{\mathcal{G}}
\newcommand{\T}{\mathcal{T}}

\newcommand{\E}{\mathbb E}

\newcommand{\ra}{\rightarrow}
\newcommand{\Lra}{\Leftrightarrow}

\newcommand{\textd}{\mathrm{d}}
\newcommand{\texte}{\mathrm{e}}

\newcommand{\cg}{c_{\mathrm{g}}}
\newcommand{\cS}{c_{\mathrm{s}}}

\newcommand\cc{\mathrm{c}}
\newcommand{\eqref}[1]{(\ref{#1})}
\makeatother

\begin{document}
\begin{frontmatter}

\title{Gibbs measures on permutations over one-dimensional discrete point sets\thanksref{T1}}
\runtitle{Gibbs measures on permutations}

\begin{aug}
\author[A]{\fnms{Marek}~\snm{Biskup}\corref{}\ead[label=e1]{biskup@math.ucla.edu}}
\and
\author[B]{\fnms{Thomas}~\snm{Richthammer}\ead[label=e2]{richth@cs.uni-hildesheim.de}}
\thankstext{T1}{Supported in part by
NSF Grant DMS-11-06850,
NSA Grant H98230-11-1-0171 and GA\v{C}R project P201-11-1558.}
\runauthor{M. Biskup and T. Richthammer}
\affiliation{UCLA and Universit\"at Hildesheim}
\address[A]{Department of Mathematics\\
UCLA\\
Los Angeles, California 90095-1555\\
USA\\
and\\
School of Economics\\
University of South Bohemia\\
Studentsk\'a 13\\37005 \v Cesk\'e Bud\v ejovice\\
Czech Republic\\
\printead{e1}} %adresu isvedimo komanda gale!
\address[B]{Institute f\"ur Mathematik\\
\quad und Angewandte Informatik\\
Universit\"at Hildesheim\\
Samelsonplatz 1\\
31141 Hildesheim\\
Germany\\
\printead{e2}}
\end{aug}

% HISTORY:
\received{\smonth{10} \syear{2013}}
\revised{\smonth{1} \syear{2014}}

% ABSTRACT
%
\begin{abstract}
We consider Gibbs distributions on permutations of a locally finite
infinite set
$X \subset\mathbb R$, where a permutation~$\sigma$ of $X$ is assigned
(formal) energy $\sum_{x \in X} V(\sigma(x)-x)$.
This is motivated by Feynman's path representation of the quantum Bose gas;
the choice $X := \mathbb Z$ and $V(x) := \alpha x^2$ is of principal interest.
Under suitable regularity conditions on the set $X$ and the potential $V$,
we establish existence and a full classification of the infinite-volume
Gibbs measures
for this problem, including a result on the number of infinite cycles
of typical permutations.
Unlike earlier results, our conclusions are not limited
to small densities and/or high temperatures.
\end{abstract}

% KEYWORDS


\end{frontmatter}

%s1 #&#
\section{Introduction}

%s1.1 #&#
\subsection{Motivation}

One of the principal difficulties underlying quantum statistical mechanics
is the noncommutative nature of the relevant observables.
For some systems, the said difficulty can sometimes be reduced
by developing a suitable classical, and often probabilistic, representation
of the problem at hand. Interestingly, this can be done
for quite a few examples of interest, namely, the quantum rotator,
the quantum Heisenberg model, the Ising model in a transversal field and,
most notably, the Bose gas.
The classical representation is still hard to analyze, but some results
often follow.
See T\'oth \cite{Toth} and Aizenman and Nachtergaele \cite{AN} for
early studies of such representations.

In this paper we take up a model that is derived
from the classical representation of interacting Bose gas.
This representation is originally due to Feynman \cite{Fe} to whom it served
as a mathematical tool for the analysis of the onset of Bose--Einstein
condensation
in ${}^4_2$He. Feynman's representation yields a model on $N$ classical
particles
in positions $x_1,\ldots,x_N \in\R^d$ that are given the weight
%
%e1.1 #&#
\begin{eqnarray}
\label{FKrep} &&\quad\sum_{\si\in\Om_X} \int W_X^{\be,\si}(
\textd B) \nonumber\\[-8pt]\\[-8pt]
&&\quad\hphantom{\sum_{\si\in\Om_X} \int}{}\times\exp \biggl\{ -\sum_{x \in X} V \bigl(\si(x)-x
\bigr) - \sum_{1 \le i < j \le N} \int_0^\be
w \bigl(B_t^{(i)} - B_t^{(j)} \bigr)
\,\textd t \biggr\}.\nonumber
\end{eqnarray}
Here $X := \{x_1,\ldots,x_N\}$ is the set of all particles, $\Om_X$ is
the set
of all permutations of $X$ (i.e., all one-to-one maps of $X$ onto $X$),
$\be$ is (twice) the inverse temperature and $V(x) := \frac{1} {2\be} |x|^2$
so that the first sum in the exponent yields independent Gaussian factors.
The probability measure $W_X^{\be,\si}$ is over collections of
independent Brownian bridges $B := (B_t^{(i)}\dvtx 0 \le t \le\be, 1
\le i \le N)$,
where the $i$th bridge starts at $x_i$ and terminates at $\si(x_i)$.
The function $w\dvtx \R^d \to\R$ is the two-body interaction potential
between the bosons. The a priori measure on the positions $x_1,\ldots,x_N$
is Lebesgue over a finite set; integrating the weight \eqref{FKrep}
defines the normalizing constant called the partition function.

A natural first case to explore is that of no interaction, that is, $w
:= 0$.
Averaging the positions over, say, a torus can then be exactly carried out
with the help of Fourier representation and the sum over permutations
can then also
be performed. An outcome of this, envisioned already by Feynman \cite{Fe},
is as follows:
For $d = 1,2$, any finite density of particles and any $\be> 0$,
a typical $\si$ will decompose entirely into finite cycles (i.e., of sizes
not growing with $N$). On the other hand,
for $d \ge3$ and each $\be> 0$ there is a critical density above which
a particle is contained in a cycle of length of order $N$ with positive
probability.
A mathematical proof of this has been given only recently
by S\"ut\H{o} \cite{S1,S2}; the critical density turns out to coincide
with that for the appearance of the Bose--Einstein condensate.

It has subsequently been observed by Betz and Ueltschi \cite{BU}
that a similar calculation to the one just mentioned can be carried out
for $V(x) = |x|^2$ replaced by more general potential functions.
The principal next challenge from the mathematics point of view is thus
to either allow for nonzero interactions, $w \neq0$, or to drop the
integration
over the positions $x_i$. The former choice is that of the prime
interest for physics; unfortunately, at this moment we do not see
any tangible way to tackle it.
The latter option is nonetheless interesting as well; it leads to
natural measures
on partitions of point sets in $\R^d$. This motivation was the basis
of an earlier article of Fichtner \cite{Fi}. More recently, Gandolfo,
Ruiz and Ueltschi
\cite{GRU} and Betz and Ueltschi \cite{BU} proposed a similar model
with the particles placed at the vertices of the integer lattice $\Z^d$.

%s1.2 #&#
\subsection{Main questions}
We will henceforth focus on the latter case and formalize it
as the following problem:
For a given locally-finite set $X \subset\R^d$ of positions and
a potential function $V\dvtx \R^d \to\R$, we wish to consider
a probability measure $\mu$ that is formally given by
%
%e1.2 #&#
\begin{equation}
\label{formalmeasure} \mu\bigl(\{\si\}\bigr) = \frac{1} Z \exp \biggl(- \sum
_{x \in X} V\bigl(\si (x)-x\bigr) \biggr),\quad\quad \si\in
\Om_X.
\end{equation}
If $X$ is finite, this corresponds to the measure arising from \eqref{FKrep}
with $w := 0$. If $X$ is infinite,
the expression \eqref{formalmeasure} is generally ill defined. In order
to extend it to infinite volume
(i.e., infinite number of particles), one either has to appeal to
limits---an approach previously used in this context by Fichtner \cite
{Fi} and Gandolfo, Ruiz and Ueltschi \cite{GRU}---or go directly by
prescribing infinite-volume Gibbs measures
via a family of specifications (Georgii \cite{G}). This will be our choice,
so the first question to ask is:
\begin{longlist}[(4)]
\item[(1)] Can we define a consistent family of specifications
of the form \eqref{formalmeasure}?
\end{longlist}
As we will see later, already this represents a departure
from the standard theory.
Naturally, one thus immediately adds:
\begin{longlist}[(4)]
\item[(2)] Under what conditions are there (infinite-volume) Gibbs measures
for the specifications defined in (1)?
\end{longlist}
In an approach via limits from finite volume, this boils down to
controlling tightness---the issue is that in the limiting measure, no
points get mapped to/from infinity.
This has, so far, only been accomplished under the assumption of
low density/high temperature (cf. conditions (V.3) or (5.8) in
Fichtner \cite{Fi})
or for interaction with a finite-range cutoff (Betz and Ueltschi
\cite{BU}, page 478).
In addition, all of this is only for the free boundary condition
(Fichtner \cite{Fi}, Theorems~2.2 and~3.1).

\begin{Note*}
 When this manuscript was very near its completion,
Betz \cite{Betz-new} posted a proof of tightness for the measures with
periodic boundary conditions over $X:=\Z^d$ in all $d\ge1$ assuming the
summability condition of the kind $\sum_{x \in\Z^d} \texte^{-
\delta
V(x)}< \infty$ for some $\delta\in(0,1)$.
\end{Note*}

Once the setting for Gibbs theory is fixed, the natural follow-up questions
concern the structure of the permutations that are typical samples from
these measures:
\begin{longlist}[(4)]
\item[(3)] Characterize the Gibbs measures that are trivial on the tail
sigma field
(i.e., the extremal Gibbbs measures).
\item[(4)] Under what conditions does $\si$ contain finite cycles only,
and when do infinite cycles occur with positive probability?
\end{longlist}
The mathematical results that are available at present all pertain to
the regime
of low densities/high temperatures: the measure defined by the free
boundary condition
contains only finite cycles almost surely (Fichtner \cite{Fi}, Theorem~3.2).
Interesting numerical simulations were performed by Gandolfo, Ruiz and Ueltschi
\cite{GRU} and Grosskinsky, Lovisolo and Ueltschi \cite{GLU}
of the model with $X := \Z^d$ (the integer lattice) and $V(x) := \al|x|^2$.
These indicated a similar dichotomy as for the ideal Bose gas:
only microscopic (finite) cycles in dimensions $d = 1,2$ and
macroscopic cycles
in $d = 3$ for $\be:= \frac{2}{\alpha}$ sufficiently large.

The principal goal of this paper is to answer the above questions
in the case of one-dimensional point sets, $X \subset\R$,
subject to natural homogeneity conditions, and a fairly rich class of
potentials $V$.
Explicitly, we show how to define Gibbs specifications,
establish the existence of a family of Gibbs measures and prove that
the extremal ones
are in one-to-one correspondence with an integer parameter called the flow.
This is the quantity whose absolute value gives the (a.s.-constant) number
of infinite cycles while its sign tells us their asymptotic direction
(it turns out that all infinite cycles necessarily ``flow'' in the same
direction).

%s1.3 #&#
\subsection{Further related work}
Apart from the above mentioned works that are focused on the connection
to the Bose gas, there are numerous other studies of combinatorial
nature that deal with similar problems. One line of research concerns
compositions of random transpositions (e.g., Schramm \cite{Schramm},
Berestycki and Durrett \cite{BerDurr}, Berestycki \cite
{Berestycki1,Berestycki2}); these pertain to situations without
underlying geometry. Another direction concerns the random stirring
process where transpositions ``arrive'' randomly but only over edges of
an underlying graph. For the graph being an infinite regular tree,
progress has recently been made concerning the existence and uniqueness
of a transition from a regime without infinite cycles to a regime with
infinite cycles (Hammond \cite{Hammond1,Hammond2}). A recent review by
Golschmidt, Ueltschi and Windridge \cite{GUW} gives further connections
between the combinatorial models and quantum systems.

There are also several alternative approaches to statistical mechanics
of Bose gases to the one proposed by Feynman. Recently, much progress
has been achieved
in the analysis of the so-called Gross--Pitaievski limit;
this density-function approach is summarized in the monograph
by Lieb, Seiringer, Solovej and Yngvason \cite{LSSY}.
Other studies were put forward that expand on the ideas of
Bogoliubov \cite{Bog};
see, for instance, a review article by Zagrebnov and Bru \cite{ZB}. The
jury is still out on which of these approaches is best suited for
understanding the physics, although some connection of the present
problem to the others has also been made; see, for example,
Ueltschi \cite{U} and the work of Adams, Bru and K\"onig \cite{ABK1,ABK2}.

%s1.4 #&#
\subsection{Outline}
The remainder of the paper is organized as follows: In the next section
we define a suitable notion of Gibbs measures for a given set of points $X$
and potential $V$ satisfying suitable assumptions. We introduce the main
concepts for our analysis and formulate a series of lemmas containing
our findings
for Gibbs measures in the given context, concluding with a summarizing theorem.
The following sections contain the proofs for these results.

%s2 #&#
\section{Definitions and results}

Here we develop the mathematical framework of our problem and give
statements of the results. Our approach is based on the theory of Gibbs
measures but, since we work in a somewhat nonstandard setting, we will
be rather pedantic
in introducing all necessary notation.

%s2.1 #&#
\subsection{Permutations on point sets and their flow}
Our aim is to construct a measure of form \eqref{formalmeasure} on
permutations on a given countably-infinite set of points $X \subset\R
$. One may want to think
of a regular point set such as the set of all integers $\Z$, but at
this point we only
assume that $X$ is:
\begin{longlist}[(ii)]
\item[(i)] locally finite (i.e., any bounded subset of $X$ is finite) and
\item[(ii)] bi-infinite (i.e., $X$ is unbounded from above and below).
\end{longlist}
%
%Note that since $X$ is a set, we exclude the possibility that points
%appear multiple times.
%For the topology on $X$ we take the standard discrete topology.
In order to be able
to identify particular points of $X$ with respect to a given position,
we use the following notation: For $a \in\R$ and $n \ge1$ let
$a_{n},a_{-n}$ be the unique points of $X$ such that
$\# ( (a,a_{n}] \cap X  ) = n = \# ( [a_{-n},a) \cap X  )$,
where $\#(A)$ is our notation for the cardinality of $A$. In other
words, $a_{n},a_{-n}$ are the $n$th point of $X$ lying (strictly) to
the right and left
of $a$, respectively. If $a \in X$ we also write $a_{0} := a$.
We set $\La^\cc:= X \setminus\La$ for $\La\subset X$,
and write $\La\Subset X$ if $\La$ is a finite subset of $X$. We will
also write
%
%e2.1 #&#
\begin{equation}
X^\star:= \biggl\{ \frac{x_{0}+x_{1}}2\dvtx  x \in X \biggr\}
\end{equation}
to denote the dual set of points of $X$.

The configuration space $\Om_X$ is the set of all permutations (i.e.,
bijections) on~$X$.
For given $\si\in\Om_X$ it will be useful to think of a pair $(x,\si(x))$
with $x \in X$ as a jump from $x$ to $\si(x)$ with the notation
%
%e2.2 #&#
\begin{eqnarray}
  x \ra y \quad&:\Lra &\quad\si(x) = y,\nonumber
\\[-8pt]\\[-8pt]
% \mbox{and}
x \leftrightarrow y \quad&:\Lra& \quad\si(x) = y\mbox{ or } \si(y) = x.\nonumber
\end{eqnarray}
We will say that $(x,\si(x))$ is a jump over $a \in\R$ if
$x < a < \si(x)$ or $\si(x) < a < x$, and it is a jump to the
right, respectively, left if $\si(x) > x$, respectively, $\si(x) < x$. The quantity
$|\si(x)-x|$ will be referred to as the length of the jump.

An important tool in our analysis will be the flow $F_a(\si)$ of a
permutation $\si$
through $a \in X^\star$. To define this object, we set
%
%e2.3 #&#
\begin{eqnarray}
  F_a^+(\si) &:=& \# \bigl\{x
\in X\dvtx  x < a < \si(x) \bigr\},\nonumber
\\[-8pt]\\[-8pt]
F_a^-(\si) &:=& \# \bigl\{x \in X\dvtx \si(x) <
a < x \bigr\}\nonumber   %
\end{eqnarray}
and define
%
%e2.4 #&#
\begin{equation}
F_a(\si) := \cases{ F_a^+(\si) - F_a^-(
\si), &\quad\mbox{if }$F_a^+(\si),F_a^-(\si) < \infty$,
\cr
\infty, &\quad\mbox{otherwise}. }
\end{equation}
As we will see, the formal (second) value is a proviso that will turn out
to be irrelevant for the typical permutations to be considered later.
A key fact is that $F_a(\si)$ does not depend on $a$:

%le2.1 #&#
\begin{Lem} \label{lemsameflow}
For locally finite $X$ and $\si\in\Om_X$,
$F_a(\si)$ has the same value for all $a \in X^\star$.
\end{Lem}

Therefore, we may (and will) drop $a$ from the notation and
define $F\dvtx \Om_X \to\Z\cup\{\infty\}$ to be the common value
of $F_a$ for all $a\in X^\star$.

Every permutation $\si\in\Om_X$ can be decomposed into disjoint
cycles, some of which
may be infinite. If $x \in X$ belongs to an infinite cycle and
$\si^n(x) \to-\infty$ for $n \to- \infty$ and  $\si^n(x) \to
\infty$
for $n \to\infty$,
we say that the cycle is going from $-\infty$ to $\infty$. Similarly we
may have cycles going
from $\infty$ to $\infty$, $-\infty$ to $-\infty$, $\infty$
to $-\infty$.
If one or both of the above limits do not exist,
we say that the infinite cycle is indeterminate.
Not too surprisingly, the value of the flow gives some information on
the number of such infinite cycles:

%le2.2 #&#
\begin{Lem} \label{lemflowcycles}
Assume that $X$ is locally finite. Any $\si\in\Om_X$ with $F(\si) =:
n \in\Z$
does not have indeterminate infinite cycles, and it has at least
$|n|$ infinite cycles from $-\infty$ to $\infty$ if $n > 0$ and at least
$|n|$ infinite cycles from $\infty$ to $-\infty$ if $n < 0$.
\end{Lem}

We note that this lemma still leaves the possibility of having
additional infinite cycles that do not
contribute to the flow, for example, an infinite cycle from $\infty$
to $\infty$ or
a pair of cycles, one from $-\infty$ to $\infty$ and one from $\infty$
to $-\infty$. These will be effectively ruled out in Theorem~\ref{thm}.

% THIS BELONGS SOMEWHERE ELSE
%Our interest in infinite cycles arises from the physics motivation for
%our investigation; indeed, the existence of infinite cycles should be
%related to the appearance of the Bose-Einstein condensate.

%s2.2 #&#
\subsection{Energy of permutations}
The energy of a permutation $\si$ will be defined in terms of a given
potential function $V\dvtx \R\to\R$. The most interesting choice is
$V(x) := \al x^2$ for $\al> 0$, but for now we will only assume $V$
to be symmetric [$V(x) = V(-x)$ for all $x \in\R$] and strictly convex.
The formal Hamiltonian
%
%e2.5 #&#
\begin{equation}
H(\si) := \sum_{x \in X} V\bigl(\si(x)-x\bigr)
\end{equation}
does not converge to a finite limit for most permutations $\si$.
As usual this problem can be avoided by considering a local version
of the energy for a given boundary configuration.
For $\La\Subset X$ we define a compatibility relation $\sim_\La$
between configurations $\si\in\Om_X$
and $\eta\in\Om_X$ by setting
%
%e2.6 #&#
\begin{equation}
\si\sim_\La\eta \quad:\Leftrightarrow\quad \forall x \in\La^\cc
\dvtx  \si(x) = \eta(x), \si^{-1}(x) = \eta^{-1}(x).
\end{equation}
In particular, if $\sigma\sim_\Lambda\eta$, then $\si$ maps $\La
\cap
\eta^{-1}(\La) = \La\cap\si^{-1}(\La)$
bijectively onto $\La\cap\eta(\La) = \La\cap\si(\La)$.
The Hamiltonian of $\si\in\Om_X$ in $\La\Subset X$ is defined by
%
%e2.7 #&#
\begin{equation}
H_{\La}(\si) := \sum_{x \in\si^{-1}(\La) \cap\La} V\bigl(\si(x)-x
\bigr).
\end{equation}
As usual, a configuration will be called a ground state of $H$
if its energy is smaller than that of any local perturbation thereof:

%de2.3 #&#
\begin{definition}\rm
$\tau\in\Om_X$ is said to be a ground state of $H$ if and only if
%
%e2.8 #&#
\begin{equation}
H_{\La}(\tau) \le H_{\La}(\si) \qquad\mbox{for all }\La\Subset X
\mbox{ and all } \si\sim_\La\tau.
\end{equation}
\end{definition}

%It turns out that our model has ground states of distinct levels
%of the specific energy: \bf SPEC ENER MAY NOT BE DEFINED \rm
It turns out that the ground states can be explicitly described:

%le2.4 #&#
\begin{Lem}
\label{lemground}
If $X$ is locally finite and bi-infinite,
and $V$ is strictly convex, then the ground states of $H$ form the
set $\{\tau_n\dvtx  n\in\Z\}$, where $\tau_n$ is the $n$-shift
permutation defined by
%
%e2.9 #&#
\begin{equation}
\tau_n(x) := x_{n} \qquad\mbox{for all } x \in X.
\end{equation}
\end{Lem}

Note that the ground states $\{\tau_n\dvtx  n \in\Z\}$ are precisely
the increasing bijections
of~$X$. [We say that $\sigma$ is strictly increasing if $x<y$ implies
$\sigma(x)<\sigma(y)$.] Moreover, they are completely parametrized by
their flow, $F(\tau_n) = n$, which (as we will see later) will be true
even for the Gibbs measures.
The proof of Lemma~\ref{lemground} and many of the following results rely
on the energy comparison of a permutation and its perturbation at
exactly two points:

%de2.5 #&#
\begin{definition}
\label{localperturb}
For $x,y \in X$ and $\si\in\Om_X$, we define $\si_{xy} \in\Om_X$
by setting
%
%e2.10 #&#
\begin{eqnarray}
  \si_{xy}(x) &:=& \si(y),\quad\quad
\si_{xy}(y) := \si(x),\nonumber
\\[-8pt]\\[-8pt]
% \mbox{
%and
%}
\si_{xy}(z) &:=& \si(z) \quad\quad\mbox{for all
} z \neq x,y. \nonumber  %
\end{eqnarray}
\end{definition}

We will sometimes refer to the transformation $\sigma\mapsto\sigma
_{xy}$ as a swap.
Note that, for $x,y,\si(x),\si(y) \in\La\Subset X$ we have $\si_{xy}
\sim_\La\si$ and
%
%e2.11 #&#
\begin{eqnarray}
&&\quad H_{\La}(\si)- H_{\La}(\si_{xy})\nonumber
\\[-8pt]\\[-8pt]
&&\quad\quad\quad= V\bigl(\si(x)-x\bigr) + V\bigl(\si(y)-y\bigr) - V\bigl(\si(x)-y\bigr) - V
\bigl(\si(y)-x\bigr).\nonumber
\end{eqnarray}
This relation is the reason why strict convexity of $V$ is crucial for
the validity of Lemma~\ref{lemground}.

%s2.3 #&#
\subsection{Specifications and Gibbs measures}
We continue to suppose that $X$ and $V$ satisfy the above assumptions.
Before defining Gibbs measures on $\Om_X$, we need to impose topological
and measurable structures on $\Om_X$.
We endow $\Om_X$ with the smallest topology under which all the projections
%
%e2.12 #&#
\begin{equation}
\si\mapsto P^+_x(\si) := \si(x) \quad\mbox{and}\quad \si\mapsto
P^-_x(\si) = \si^{-1}(x) \quad\quad(x \in X)
\end{equation}
are continuous.
[Identifying $\si$ with $(\si(x),\si^{-1}(x))_{x \in X}$ this
topolgy on
$\Om_X$ coincides with the product of the   discrete topologies on $X
\times X$.]
This topology is metrizable, for example, by
%
%e2.13 #&#
\begin{equation}
d(\si,\eta) := \inf \bigl\{2^{-r}\dvtx \si= \eta\mbox{ and }
\si^{-1} = \eta^{-1} \mbox{ on } X \cap[-r,r] \bigr\},
\end{equation}
and $\Om_X$ is thus a complete separable metric space.

Let $\F_X := \si(P^+_x,P^-_x\dvtx  x \in X)$ denote the Borel-$\si
$-algebra on $\Om$,
and for $\La\subset X$ we use
%
%e2.14 #&#
\begin{equation}
\F_{\La} := \si\bigl(P^+_x,P^-_x\dvtx  x \in
\La\bigr)
\end{equation}
to denote the $\si$-algebra of events depending on $\La$ only.
A function $f\dvtx \Om_X \to\R$ is called local if it is measurable
with respect to  $\F_{\La}$
for some $\La\Subset X$; an event is called local if its indicator is
a local function. Since every local event is open and closed,
every local function is continuous. As usual we can use these to define
%
%e2.15 #&#
\begin{equation}
\T:= \bigcap_{\La\Subset X} \F_{\La^\cc},
\end{equation}
the tail-$\si$-algebra of all events that do not depend on what a
permutation looks like on any bounded set.

In order to construct (infinite volume) Gibbs measures for the
above Hamiltonian $H$, we will use the method of specifications.
Here, the specification corresponding to boundary configuration
$\eta\in\Om_X$ and volume $\La\Subset X$ is the discrete probability
measure $\ga_{\La}(\cdot|\eta)$ on $(\Om_X,\F_X)$ defined by
%
%e2.16 #&#
\begin{eqnarray}
\ga_{\La}\bigl(\{\si\}|\eta\bigr) &:=&
\frac{1}{Z_{\La}(\eta)} \texte^{- H_{\La}(\si)}1_{\{\si\sim
_\La
\eta\}},\nonumber
\\[-8pt]\\[-8pt]
\eqntext{\mbox{where } Z_{\La}(\eta) := \sum_{\si\colon\si\sim_\La\eta}
\texte^{-
H_{\La
}(\si)}.}
\end{eqnarray}
We note that  $Z_{\La}(\eta) >0$ (since $\si:= \eta$ gives a nonzero
contribution)
and $Z_{\La}(\eta) < \infty$ (since only finitely many permutations
are $\sim_\La$-compatible to $\eta$ if $\La$ is finite).
%WHY IS THIS IMPORTANT?
%Thus $\ga_{\La}$ is a well defined probability kernel from $(\Om_X,\F_{
Obviously, $\eta\mapsto\gamma_\Lambda(A|\eta)$ is measurable for
each $A\in\F_X$.

%de2.6 #&#
\begin{definition}
A Gibbs measure $\mu$ with respect to the above Hamiltonian $H$ is a
probability
measure on $(\Om_X,\F_X)$ that is compatible with all specifications:
%
%e2.17 #&#
\begin{equation}
E_\mu \bigl(\ga_{\La} (A|\cdot) \bigr) = \mu(A) \quad\quad\mbox{for
all $A \in\F_X$ and all $\La\Subset X$},
\end{equation}
where $E_\mu$ denotes expectation with respect to $\mu$.
We will write $\G$ to denote the set of all Gibbs measures for the
(implicit) Hamiltonian $H$.
\end{definition}

%re2.7 #&#
\begin{remark}
Our setup differs from that in Georgii \cite{G}, which is versed in
terms of spin systems. Such a description is permissible in our case as
well. A natural attempt would be to proclaim $\sigma(x)$ to be
an $X$-valued spin at $x$, but that choice does not lead to quasilocal
specifications. Taking a pair of values $(\sigma(x),\sigma^{-1}(x))$
for the spin at $x$ solves the quasilocality problem, but only at the
cost of introducing hard-core restrictions.
\end{remark}

%s2.4 #&#
\subsection{Conditions on the point set and the potential}
In order to show the existence (and discuss further properties) of
Gibbs measures, we have to impose conditions on $X$ and $V$
that are stronger than those considered so far.
We did not aim at finding the most general conditions under which our
results hold,
but rather gave conditions under which the proofs remain transparent
and which still include the main examples of interest.

For the
formulation of the conditions on $X$ we introduce constants
that control the spacings between the points of
$\La(a,n) := \{a_{-|n|-1},\ldots,a_{|n|+1}\}$, $(n \in\Z)$
and the growth rate of the number of particles of $X$ as seen from $a$:
%
%e2.18 #&#
\begin{equation}
\cS^X(a,n) := \min \left\{c \ge1\dvtx  %
\begin{aligned}&\mbox{ consecutive points of $\La(a,n)$}
\\[-4pt]
 &\mbox{ keep distances $\in\bigl[c^{-1},c\bigr]$} \end{aligned}
\right\}
\end{equation}
and
%
%e2.19 #&#
\begin{equation}
\qquad \cg^X(a) := \inf \bigl\{ c \ge0\dvtx \#\{x \in X\dvtx 0
< |x-a| \le t\} \le c t \mbox{ for all } t > 0 \bigr\}.
\end{equation}
(The subindices ``s,'' resp., ``g'' stand for ``separation,'' resp.,
``growth.'')
Below we will consider the following, progressively restrictive,
properties of $X \subset\R$ (and $n \ge0$):
\begin{longlist}[(X$n$)]
\item[(X1)] $X$ is locally finite and bi-infinite (just as before).
\item[(X2)] $\cg^X(a) < \infty$ for some $a \in\R$.
\item[(X$n$)] For some $c_n < \infty$, bi-infinitely many points
$a \in X^\star$ satisfy
%
%e2.20 #&#
\begin{equation}
\cS^X(a,n) \le c_n \quad\mbox{and}\quad \cg^X(a) \le
c_n.
\end{equation}
\end{longlist}
We note that $X := \Z$ satisfies the given conditions
[since $c^\Z_s(a,n)=1$ and $c^\Z_g(a) = 2$ for all $a \in X^\star$],
but the conditions also allow for point sets such that the distance
of consecutive points is not bounded from above or below.
In paricular, we can consider point sets produced by a Poisson point process
or other rather general shift-invariant processes:

%le2.8 #&#
\begin{Lem} \label{lempoisson}
Consider a point process (i.e., an $\N\cup\{\infty\}$-valued random
purely-atomic Borel measure) $\mathscr X$ on $\R$ that has the
following properties:
\begin{longlist}[(3)]
\item[(1)] $\mathscr X(\R)=\infty$ a.s. but $E\mathscr X(A)<\infty$
for any compact $A\subset\R$.
\item[(2)] $\mathscr X$ is simple, that is, a.s. no two points
of $\mathscr X$ coincide.
\item[(3)] The law of $\mathscr X$ is invariant and ergodic under the
map $x\mapsto x+1$.
\end{longlist}
Then the set of points corresponding to a.e. sample of $\mathscr X$
satisfies \textup{(X1)} and \textup{(X}$n$\textup{)} for all $n \in\Z$.
\end{Lem}

Concerning the potential function $V\dvtx \R\to\R$, we will consider
the following properties:
\begin{longlist}[(V2)]
\item[(V1)] $V$ is symmetric and strictly convex (just as before).
\item[(V2)] $V$ satisfies the following growth condition: For all $d >
0$ we have
%
%e2.21 #&#
\begin{equation}
\psi_d(x) := \frac{V(x)+V(0) - 2 V( (x+d)/2)}{x\log x} \to\infty \quad\quad
\mbox{for } x \to\infty.
\end{equation}
\end{longlist}
We note that $V(x) := \al|x|^{1+\ep}$ satisfies both (V1) and (V2)
for all $\al> 0$ and all $\ep> 0$. In particular, this includes
the most interesting case $V(x) := \al x^2$ for $\al> 0$.
Linearly growing potentials, for example, $V(x) := \al|x|$, satisfy neither
(V1) nor (V2), and indeed we are unable to extend our conclusions to
these cases.

%s2.5 #&#
\subsection{Existence of Gibbs measures}
We note that without any assumptions on $V$ it may be the case
that there are no Gibbs measures. This can already be seen when $V := 0$,
which corresponds to the physically interesting case of zero
temperature (i.e., $\alpha:=0$).

%le2.9 #&#
\begin{Lem} \label{lemnonexistence}
For infinite $X$ and $V := 0$ we have $\G= \varnothing$.
\end{Lem}

On the other hand, existence can be shown under the conditions
introduced in the
previous section. At this point the weaker condition (X2) is sufficient.

%le2.10 #&#
\begin{Lem} \label{lemexistence}
Let $X$ and $V$ satisfy \textup{(X1)}, \textup{(X2)} and \textup{(V1)},
\textup{(V2)}, respectively. For
every $n \in\Z$
there is $\mu\in\G$ such that $\mu(F = n) = 1$.
\end{Lem}

So there are in fact many Gibbs measures, at least one for each value
of the flow.
In all cases the flow is finite and this turns out to be no accident:

%le2.11 #&#
\begin{Lem}\label{lemfiniteflow}
Let $X$ and $V$ satisfy \textup{(X1)}, \textup{(X2)} and \textup{(V1)},
\textup{(V2)}, respectively.
Then $F$ is $\T$-measurable and finite $\mu$-a.s. for every $\mu\in
\G$.
\end{Lem}

As in the classical Gibbs measure theory, $\G$ is a closed convex set.
Its extremal elements are called extremal Gibbs measures, and they are
precisely those that trivialize on the tail sigma algebra $\T$.
Consequently, every Gibbs measure
can be decomposed into extremal Gibbs measures by conditioning on $\T$.
The previous lemma implies that $F$ is a.s. constant with respect to
every extremal
Gibbs measure. We may thus filter the Gibbs measures according to the
flow and focus attention on sets
%
%e2.22 #&#
\begin{equation}
\G_n := \bigl\{\mu\in\G\dvtx \mu(F=n) = 1 \bigr\}.
\end{equation}
A key challenge now is to relate the value of the flow to the number of
infinite cycles.

%s2.6 #&#
\subsection{Infinite cycles and classification of Gibbs measures}
Our analysis of infinite cycles requires introduction of an additional
technical tool, called a cut. This is defined as follows:

%de2.12 #&#
\begin{definition}
\label{def1}
Let $\si\in\Om_X$ be a permutation with flow $F(\si) =: n \ge0$
and $a \in X^\star$. Then $a$ is called a cut for $\si$ if $F_a^+(\si
)=n$ [and so $F_a^-(\sigma)=0$] and $\si$ contains the following $n$ jumps
%
%e2.23 #&#
\begin{equation}
a_{-n} \to a_{1}, \ldots, a_{-1} \to
a_{n}.
\end{equation}
For $n < 0$ similar notions are defined by reversing the directions of
the jumps.
\end{definition}

The cuts are helpful for the following reason:
If $a$ is a cut for $\si$, then $\si$ cannot have cycles jumping over $a$
apart from the $|n|$ infinite ones described in Lemma~\ref
{lemflowcycles}. In particular, we can make the following conclusions:

%
%le2.13 #&#
\begin{Lem} \label{leminfcycles}
Let $X$, $V$ and $n\in\Z$ satisfy \textup{(X1)}, \textup{(X}$n$\textup{)}
and \textup{(V1)}, \textup{(V2)}, and let
$\mu\in\G_n$. Then:
\begin{longlist}[(b)]
\item[(a)] there are bi-infinitely many cuts for $\mu$-a.e.
permutation, and
\item[(b)] $\mu$-a.e. permutation has exactly $|n|$ infinite cycles.
\end{longlist}
\end{Lem}

A duplication argument yields a uniqueness result.
For this we call $a \in X^\star$ a cut for a pair of permutations
$(\si,\si')$ if it is a cut for both $\si$ and $\si'$.

%le2.14 #&#
\begin{Lem} \label{lemuniqueness}
Let $X$, $V$ and $n\in\Z$ satisfy \textup{(X1)}, \textup{(X}$n$\textup{)}
and \textup{(V1)}, \textup{(V2)},
and let $\mu,\mu' \in\G_n$. Then:
\begin{longlist}[(b)]
\item[(a)] there are bi-infinitely many cuts for
$\mu\otimes\mu'$-almost every pair of permutations, and
\item[(b)] $\mu= \mu'$.
\end{longlist}
\end{Lem}

We conclude this section by collecting all of our previous findings and
thus providing
a complete description of the set of all Gibbs measures.

%th2.15 #&#
\begin{Thm}
\label{thm}
Let $X$, $V$ and $n\in\Z$ satisfy \textup{(X1)}, \textup{(X}$n$\textup{)}
and \textup{(V1)}, \textup{(V2)}. Then the
following holds:
\begin{longlist}[(b)]
\item[(a)] $\G_n$ contains a single Gibbs measure $\mu_n$;
\item[(b)] $\mu_n$-a.e. permutation has exactly $|n|$ infinite cycles;
\item[(c)] for any $\eta\in\Om$ with $F(\eta) = n$ and any sequence
of increasing subsets $\La\uparrow X$, $\mu_n$ is the weak limit of
specifications $\ga_\La(\cdot|\eta)$.
\end{longlist}
In particular, if \textup{(X}$n$\textup{)} is satisfied for every
$n \in\Z$ [and
\textup{(V1)}, \textup{(V2)} hold], then
%
%e2.24 #&#
\begin{equation}
\G= \biggl\{\sum_{n \in\Z} c_n
\mu_n \dvtx  c_n \ge0, \sum_n
c_n = 1 \biggr\}.
\end{equation}
\end{Thm}

We finish with a few remarks:

%re2.16 #&#
\begin{remark}
As already mentioned, our assumptions on $V$ are by far not optimal. In
fact, we believe that a majority of our results---that is, with the
exception of the characterization of ground states in Lemma~\ref
{lemground}---carry over when $V\dvtx \R\to\R$ is perturbed by
adding a
continuous, even function that decays sufficiently fast to zero at infinity.
\end{remark}

%re2.17 #&#
\begin{remark}
Our main theorem states that the specifications $\gamma_{\La}(\cdot
|\eta
)$ have a weak limit along volumes increasing to $X$ for any boundary
condition $\eta$ that has a finite flow.
If $\eta$ has infinite flow, this is not necessarily the case.
To illustrate different behaviors we give two examples, assuming that
(V1) and (V2) are satisfied and that $X := \Z$ for simplicity.
\end{remark}

%ex1 #&#
\begin{Example}\label{ex1}
Let $\eta\in\Om_\Z$ be defined by $\eta(0) := 0$ and $\eta(x) :=
x +
2p_x$, where $p_x$ is the maximal power of 2 dividing $|x|$. Writing
$x:=(2n+1)2^k$ for $x\ne0$ and some $n\in\Z$ and $k\in\N$, we get
$\eta
(x)=(2n+3)2^k$ and so $\eta$ is a permutation. It is easy to check
that $\eta$ has infinitely many cycles from $-\infty$ to $\infty$.
Here an argument underlying the proof of Lemma~\ref{lemnestedjump}
can be used to show that for all $x,y \in\Z$ we have $\ga_\La(x \to
y|\eta) \to0$
whenever $\La\uparrow\Z$. (Indeed, with increasing $\Lambda$ there is
an increasing number of jumps over $[x,y]$ with both endpoints
in $\Lambda$; any of these can be swapped with the jump $x\to y$ while
gaining energy.)
So in this case no sequence of specifications (with $\Lambda\uparrow
\Z
$) is tight and no weak limits can be extracted.
\end{Example}

%ex2 #&#
\begin{Example}\label{ex2}
Let $\eta\in\Om_\Z$ be defined by $\eta(x) := -x$
for all $x \in\Z$. Then $\eta$ has infinite flow, but no infinite cycles.
Since the restriction of $\ga_\La(\cdot|\eta)$ to $\La\cap(-\La
)$ is
the same
as $\ga_{\La\cap(-\La)}(\cdot|\tau_0)$, here we have $\ga_\La
(\cdot
|\eta) \to\mu_0$ weakly whenever $\La\uparrow X$.
\end{Example}

%re2.18 #&#
\begin{remark}
The one-dimensional nature of the underlying set $X$ has been crucial
for our reasoning, mainly due to the concept of the flow. Utilizing
this concept
we are able to show that for any given jump of sufficient length
it is possible to find another jump such that the corresponding swap
(see Definition~\ref{localperturb}) decreases the energy significantly.
In two dimensions it is not clear how to define a useful quantity analogous
to the flow, and indeed the idea of reducing the energy by a suitable
swap fails: If a long jump is entirely surrounded by jumps of the same length
and in the same direction, swapping this jump with another one may
in fact increase the energy. As a consequence our technique does not
even give
existence of Gibbs measures in two or more dimensions.
\end{remark}

%s3 #&#
\section{Proofs: Preliminary observations}

We are now ready to commence the exposition of the proofs. Here we
begin with some preliminary observations. Certain technical aspects
that feed into the main arguments are then discussed in Section~\ref
{sec4}; the main results are proved in Sections \ref{sec5} and \ref
{sec6}. This section gives the proofs of Lemmas \ref{lemsameflow},
\ref{lemflowcycles}, \ref{lemground} and \ref{lempoisson}.

%s3.1 #&#
\subsection{Permutations and their flow}
We begin by giving the proofs of the lemmas dealing with existential
facts and properties of the flow of a permutation.

\begin{pf*}{Proof of Lemma~\ref{lemsameflow}}
Let $\si\in\Om_X$.
It suffices to show that $F_a(\si) = F_b(\si)$ whenever $a < b \in
X^\star$
with only one point $x \in X$ in between.
Then every jump except the one to and the one from $x$ contributes to
$F_a^+(\si)$ and $F_b^+(\si)$ exactly the same way [and likewise for
$F_a^-(\si)$ and $F_b^-(\si)$]. This already implies that $F_a(\si) =
\infty$ if and only if
$F_b(\si) = \infty$.
The jumps to and from $x$ decompose into five cases:
(1)~$\si(x) = x$, (2)~$\si^{-1}(x) < x < \si(x)$, (3)~$\si(x) < x <
\si
^{-1}(x)$,
(4)~$x > \si^{-1}(x), \si(x)$ and (5)~$x < \si(x),\si^{-1}(x)$. In each
of these it is easy
to check that the contribution of the combined jumps to $F_a(\si)$ is
the same
as the one to $F_b(\si)$. Thus $F_a(\si) = F_b(\si)$ as desired.
\end{pf*}

\begin{pf*}{Proof of Lemma~\ref{lemflowcycles}}
Pick $\si\in\Om_X$ with $n:=F(\si)$ finite, and let $x\in X$ be such
that $\{\si^k(x)\dvtx  k \in\Z\}$ is an infinite cycle. Since $X$ is
locally finite, and $\si^k(x)$ are all distinct, $\sigma^k(x)$
eventually leaves every bounded interval for both $k \to\pm\infty$.
Moreover, as $F_a(\si) = n$ for any $a \in X^\star$, the
sequence $\sigma^k(x)$ can jump over $a$ at most finitely many times
and $\sigma^k(x)-a$ changes sign only finitely often. Hence, the limits
of $\sigma^k(x)$ as $k\to\pm\infty$ exist in $\{\pm\infty\}$. This
rules out indeterminate infinite cycles.

In order to compare the number of infinite cycles of $\si$ with $n =
F_a(\si)$,
consider all cycles of $\si$ and note that finite cycles, cycles
from $\infty$ to $\infty$
and cycles from $-\infty$ to $-\infty$ do not contribute to $F_a(\si
)$ since
they jump over $a$ equally (finitely) often to the left and to the right.
Cycles from $-\infty$ to $\infty$ contribute $1$ to $F_a(\si)$ since
they jump
to the right one more times than to the left. Similarly,
cycles from $\infty$ to $-\infty$ contribute $-1$ to $F_a(\si)$.
It follows that there have to be at least $|n|$ infinite cycles
of the same orientation as the sign of $n$.
\end{pf*}

%s3.2 #&#
\subsection{Ground states}
Using the local perturbation from Definition~\ref{localperturb}
we easily get the following monotonicity result for local minima of the energy:

%le3.1 #&#
\begin{Lem} \label{lemminenergy}
Let $X$ be locally finite and $V$ be strictly convex.
Let $\La\Subset X$ and $\si\in\Om$. The minimum
$\min\{H_{\La}(\tau)\dvtx \tau\sim_\La\si\}$ is attained at the unique
$\tau\sim_\La\si$ such that $\tau\dvtx \La\cap\si^{-1}(\La)
\to
\La\cap\si(\La)$
is a (strictly) increasing function.
\end{Lem}

\begin{pf}
There are only a finite number of permutations $\tau$ with $\tau\sim
_\La\si$ and so the minimum of $H_{\La}$ is attained.
Moreover, any $\tau\sim_\La\si$ maps $\La\cap\si^{-1}(\La) =
\La
\cap\tau^{-1}(\La)$
bijectively onto $\La\cap\si(\La) = \La\cap\tau(\La)$ and it
coincides with
$\si$ elsewhere. If $\tau\dvtx \La\cap\si^{-1}(\La) \to\La\cap
\si
(\La)$
is not an increasing function, there are $x,y \in\La\cap\tau
^{-1}(\La)$
such that $x < y$ and $\tau(x) > \tau(y)$. As $x,y,\tau(x),\tau(y)
\in
\La$, we have $\tau_{xy} \sim_\La\tau\sim_\La\si$. We also note that
%
%e3.1 #&#
\begin{eqnarray}
&&\quad H_{\La}(\tau) -H_{\La}(\tau_{xy})\nonumber
\\[-8pt]\\[-8pt]
&&\quad\quad\quad= V\bigl(\tau(x)-x\bigr) + V\bigl(\tau(y)-y\bigr) - V\bigl(\tau(x)-y\bigr) - V
\bigl(\tau(y)-x\bigr) > 0,\nonumber
\end{eqnarray}
where the inequality follows from the strict convexity of $V$ and the
assumed ordering of $x$, $y$, $\tau(x)$ and $\tau(y)$; see Lemma~\ref
{lemma-4.2}.
Thus $H_{\La}$ does not take its minimum at a nonincreasing $\tau$.
The minimizing $\tau$ is unique because it is a bijection of $\La\cap
\si^{-1}(\La)$ onto $\La\cap\si(\La)$, and is determined
by $\sigma$
elsewhere.
\end{pf}

\begin{pf*}{Proof of Lemma~\ref{lemground}}
First, $\tau_n$ is a ground state by Lemma~\ref{lemminenergy},
since $\tau_n$ is increasing on $\Lambda\cap\tau_n^{-1}(\Lambda)$ for
every $\La\Subset X$.
On the other hand, if $\tau\neq\tau_n$ for all $n$, then $\tau$ is
not increasing;
that is, there are $x,y$ such that $x<y$ and $\tau(x) > \tau(y)$.
Considering $\La\Subset X$ such that $x,\tau(x),y,\tau(y) \in\La$,
Lemma~\ref{lemminenergy} shows that $H_\Lambda(\tau)$ is not minimal
and so $\tau$ is not a ground state.
\end{pf*}

%s3.3 #&#
\subsection{Point sets arising from point processes}
Here we show that samples from point processes fulfilling the premises
of Lemma~\ref{lempoisson} automatically satisfy the requirements (X1)
and (X$n$) for any $n\in\Z$.
Recall that a point process $\mathscr X$ on $\R$ is a random
purely-atomic $\N\cup\{\infty\}$-valued Borel measure on $\R$. Assuming
that compact sets receive finite mass almost surely, $\mathscr X$ can
be interpreted as a sum of unit point masses. If, furthermore, no
points are degenerate a.s., which technically means that $\mathscr X(\{
x\})\in\{0,1\}$ for every $x$, then samples from $\mathscr X$ can be
identified with point sets.

\begin{pf*}{Proof of   Lemma~\ref{lempoisson}}
Let us first prove (X1). Condition (1) ensu\-res that
$\mathscr X$ is locally finite a.s. Since also $\mathscr X(\R)=\infty$
a.s., $\mathscr X$ contains infinitely many
points a.s. We cannot have $\mathscr X([0,\infty))<\infty$ with a
positive probability, because the shift invariance
of the law of $\mathscr X$ would also imply $\mathscr X([-L,\infty
))\stackrel{\mathrm{law}}{=}\mathscr X([0,\infty))<\infty$ with the
same probability for all $L\in\N$ and, taking $L\to\infty$,
also $\mathscr X(\R)<\infty$ with the same probability, in
contradiction with (1). Similarly $\mathscr X((-\infty,0]) = \infty$,
and we conclude that (X1) holds a.s.

Moving over to (X$n$), we fix $n \in\N$, pick constants $K_1,K_2$ and
introduce
the events
%
%e3.2 #&#
\begin{equation}
A_k^{(1)}:= \biggl\{\sup_{\ell\ge1}
\frac{\mathscr X([k-\ell,k+\ell
))}\ell \le K_1 \biggr\} \quad\mbox{and}\quad
A_k^{(2)}:= \bigl\{ c_s^{\mathscr X} (k,n)
\le K_2 \bigr\}\hspace*{-30pt}
\end{equation}
for $k \in\Z$.
We claim that $K_1,K_2$ can be chosen so large that each of these
events has probability in excess of
$\frac{1} 2$. For the first we observe that\break $\frac{\mathscr X([k-\ell
,k+\ell))}{2\ell} \to
\E\mathscr X([k,k+1))$ a.s. by the Birkhoff pointwise ergodic theorem
(using the additivity of $\mathscr X$ and shift-ergodicity of the law
of $\mathscr X$) and the limit is finite by
assumption~(1). For the second we note that a.s. there are $n+1$
points of $\mathscr X$ to the right
of $k$ and $n+1$ points to the left of $k$ since $\mathscr X$ is
bi-infinite, and together with
the nondegeneracy of points from condition (2) this implies that the
distances between
the points of $\La(k,n)$ are positive and finite.

Now set $A_k:=A_k^{(1)}\cap A_k^{(2)}$ and note that $A_k$ has positive
probability. So, by ergodicity and the Birkhoff Theorem again, $A_k$
occurs at bi-infinitely many---in fact, a positive density of---$k$'s
almost surely. Suppose $A_k$ occurs and set $a_k := \frac{k_1
+k_{-1}}2$ if $k \notin\mathscr X$ and $a_k := \frac{k_1 +k_0}2$ if $k
\in\mathscr X$ so that $a_k \in\mathscr X^\star$.
We observe that $c_s^{\mathscr X}(a_k,n) \le K_2$ by containment in $A_k$,
so it suffices to bound $c_g^{\mathscr X}(a_k)$.
Since $|k-a_k| \le\frac{K_2} 2$ we have
%
%e3.3 #&#
\begin{eqnarray}
\label{boundset} \#\bigl\{x \in\mathscr X \dvtx
|x-a_k| \le t\bigr\}&\le& \#\biggl
\{x \in\mathscr X \dvtx |x-k|\le  t +   \frac{K_2} 2 \biggr\} \nonumber\\[-8pt]\\[-8pt] &\le&
\biggl( t +  \frac{K_2} 2  +1\biggr) K_1\nonumber
\end{eqnarray}
using containment in $A_k$. Thus for $t\ge\frac{1} {2K_2}$ we can estimate
%
%e3.4 #&#
\begin{eqnarray}
  \frac{1}t \#\bigl\{x \in\mathscr X\dvtx 0 <
|x-a| \le t\bigr\} &\le& \frac{t +   K_2/2  +1}{t} K_1\nonumber
\\[-8pt]\\[-8pt]
&\le& 2K_2 \biggl( \frac{1} {2K_2}  +
\frac{K_2} 2 +1  \biggr) K_1 \nonumber  %
\end{eqnarray}
using monotonicity. For $t < \frac{1} {2K_2}$ the expression on the left
is $0$ by the
lower bound on the gap between successive points.
This shows $\cg^{\mathscr X}(a_k)\le2K_2 ( \frac{1} {2K_2}  +
  \frac{K_2} 2+1 ) K_1$.
These bounds apply at bi-infinitely many $a\in\mathscr X^\star$, so
the condition (X$n$) is successfully verified.
\end{pf*}

%s4 #&#
\section{Technical issues}
\label{sec4}
In this section we collect some straightforward consequences
of our assumptions on $X$ and $V$. Then we proceed to discuss how to
obtain energy estimates and how these can be converted to bounds on
probabilities of events.

%s4.1 #&#
\subsection{Consequences of our assumptions}
As the set $X$ will be always clear from the context, we will
henceforth write $\cS(a,n)$ for $\cS^X(a,n)$. We note that, for any
discrete bi-infinite $X \subset\R$, automatically $\cS(a,n) < \infty$,
since this quantity is determined by only finitely many nonzero and
finite distances.
We also observe that (X1) and (X2) imply that $\cg(a) < \infty$ for
every $a \in\R$.

%le4.1 #&#
\begin{Lem}
If $\cg(a)<\infty$ and $f\dvtx [0,\infty) \to[0,\infty)$ is
nonincreasing, then
%
%e4.1 #&#
\begin{eqnarray}
\label{sumscga} %
  \sum_{x \in X\dvtx  x > a}
f(x-a) &\le&\cg(a) \int_0^\infty\,\textd s\, f(s),\nonumber
\\[-8pt]\\[-8pt]
\sum_{x \in X\dvtx  x < a} f(a-x) &\le&\cg(a) \int
_0^\infty\,\textd s\, f(s) \nonumber  %
\end{eqnarray}
and
%
%e4.2 #&#
\begin{equation}
\label{sumscga2} \sum_{x,y \in X\colon x < a < y} f(y-x) \le
\cg(a)^2 \int_0^\infty \textd s \,s f(s).
\end{equation}
\end{Lem}

\begin{pf}
For the first estimate we note $n \le\cg(a) (a_{n} -a)$ by definition
of $\cg(a)$, so
%
%e4.3 #&#
\begin{eqnarray}
  \sum_{x \in X\colon x > a} f(x-a) &=&
\sum_{n > 0} f(a_{n}-a)\nonumber
\\
&\le&\sum_{n > 0} f \biggl(\frac{n}{\cg(a)} \biggr)
\\
&\le&\int_0^{\infty} \textd s\, f \biggl(
\frac{s}{\cg(a)} \biggr) = \cg(a) \int_0^{\infty}
\textd x\, f(x) \nonumber  %
\end{eqnarray}
using the assumed monotonicity of $f$, integral comparison for monotone
sums and the substitution
$x := s/{\cg(a)}$. The second estimate is obtained similarly; the third
can be obtained
by first using the first estimate for the sum over $y$ and then the
second estimate for the sum over $x$ to get
%
%e4.4 #&#
\begin{equation}
\sum_{x,y \in X\colon x < a < y} f(y-x) \le
\cg(a)^2 \int_0^\infty \textd x\int
_0^\infty \textd y \,f(x+y).
\end{equation}
The result follows by substituting $s := x+y$.
\end{pf}

Now we turn to properties of the potential function $V$:

%le4.2 #&#
\begin{Lem}
\label{lemma-4.2}
Suppose $V\dvtx \R\to\R$ is strictly convex and symmetric. Then
%
%e4.5 #&#
\begin{equation}
\label{Vmonotone}  V \mbox{ is strictly increasing in $[0,\infty)$.}
\end{equation}
In addition, for all $t_1 < s_1 \le s_2 < t_2 \in\R$
such that $s_1+s_2 = t_1 +t_2$ we have
%
%e4.6 #&#
\begin{eqnarray}
  \label{convrev}&& V(t_1) +
V(t_2) -  V(s_1) - V(s_2)\nonumber
\\[-8pt]\\[-8pt]
 &&\quad\quad\ge2 \frac{t_2-s_2}{t_2-t_1} \biggl(V(t_2)+V(t_1) - 2 V
\biggl(\frac{t_1+t_2}{2} \biggr) \biggr)  \nonumber %
\end{eqnarray}
and thus, in particular,
%
%e4.7 #&#
\begin{equation}
\label{convst} V(t_1) + V(t_2) - V(s_1) -
V(s_2) > 0.
\end{equation}
\end{Lem}

\begin{pf}
The first conclusion is trivial. For \eqref{convrev}, comparing the
slopes of secant lines yields
%
%e4.8 #&#
\begin{eqnarray}
  \frac{V(s_1) - V(t_1)}{s_1-t_1} &\le&\frac{V((t_1+t_2)/2)-
  V(t_1)}{(t_2-t_1)/2}\nonumber
\\[-8pt]\\[-8pt]
&\le&\frac{V(t_2)-V((t_1+t_2)/2)}{(t_2-t_1)/2} \le\frac{V(t_2)
 - V(s_2)}{t_2-s_2}.  \nonumber %
\end{eqnarray}
The bound follows by estimating the difference of the outer terms
against the difference
of the inner terms and invoking $t_2-s_2 = s_1 - t_1$.
\end{pf}

If $V$ satisfies (V2), it will be convenient to denote
$c_\psi(d,N):=\sup\{x\in\R\dvtx \psi_d(x)<N\}$. Then
%
%e4.9 #&#
\begin{equation}
\label{cpsi} \psi_d(x) \ge N\quad\quad \mbox{for all } x \ge
c_\psi(d,N).
\end{equation}
This relation will be quite useful in what follows.

%s4.2 #&#
\subsection{Energy estimates}
Many of the following results rely on energy estimates, comparing
a permutation $\si$ with jumps $v \to w, y \to z$ to the swapped
configuration $\si_{vy}$. The following estimate will be helpful:

%le4.3 #&#
\begin{Lem} \label{lemenergyestimate}
Suppose $V\dvtx \R\to\R$ is strictly convex and even.
Let $\La\Subset X$ and $v,w,y,z \in\La$ and suppose that $y',z' \in
\R
$ are such
that  $v < y' \le z' < w$,  $y' \le y$ and $z \le z'$.
If $\si\in\Om_X$ contains jumps $v \to w$, $y \to z$, then
%
%e4.10 #&#
\begin{equation}
\label{energyestimate} \quad\quad H_{\La}(\si)- H_{\La}(\si_{vy}) \ge2
\min\bigl\{y'-v,w-z'\bigr\} \psi _{z'-y'}(w-v)
\log(w-v).
\end{equation}
\end{Lem}

\begin{pf}
Let $\De H := H_{\La}(\si)- H_{\La}(\si_{vy})$ denote the energy
difference to be estimated. Then
%
%e4.11 #&#
\begin{eqnarray}
  \De H &=& V(w-v) + V(z-y) - V(z-v) - V(w-y)\nonumber
\\
&\ge& V(w-v) + V\bigl(z-y'\bigr) - V(z-v) - V\bigl(w-y'
\bigr)
\\
&\ge& V(w-v) + V\bigl(z'-y'\bigr) - V
\bigl(z'-v\bigr) - V\bigl(w-y'\bigr),  \nonumber
\end{eqnarray}
where for the inequalities we used \eqref{convst} along with $w-y,z-y'
\in[z-y,w-y']$ and
$z-v,z'-y' \in[z-y',z'-v]$.
Since $z'-v,w-y' \in[z'-y',w-v]$ we now use \eqref{convrev} to obtain
%
%e4.12 #&#
\begin{eqnarray}
  \De H &\ge&2 \frac{\min\{y'-v,w-z'\}}{w-v-(z'-y')}\nonumber
\\[-8pt]\\[-8pt]
&&{}\times \biggl(V(w-v)+V\bigl(z'-y'\bigr)-2V \biggl(
\frac{w-v+z'-y'}{2} \biggr) \biggr).\nonumber   %
\end{eqnarray}
Since the terms in the large parentheses nonnegative, we now estimate
%
%e4.13 #&#
\begin{equation}
\frac{\min\{y'-v,w-z'\}}{w-v-(z'-y')} \ge\frac{\min\{y'-v,w-z'\}}{w-v}
\end{equation}
and then apply $V(z'-y') \ge V(0)$, as implied by \eqref{Vmonotone}, to
obtain \eqref{energyestimate}.
\end{pf}

If we only have control over one jump $v \to w$ of a permutation $\si$,
we need to find a suitable second jump $y \to z$ to make the above
energy estimate work. This is achieved by the following:

%le4.4 #&#
\begin{Lem} \label{lemfindjump}
Let $a \in\R$ and let $\si\in\Om_X$ with $|F(\si)| =: n \in\N$
contain a jump $(v,w)$
over $a_{1},\ldots,a_{n}$ to the right. Then $\si$ contains a jump $(y,z)$
such that $y \ge a_{1}$ and $z \le a_{n}$.
\end{Lem}

\begin{pf}
We consider the $n$ jumps starting from $a_{1},\ldots,a_{n}$. If one of
them does
not jump over $b:= \frac{1} 2 (a_{n}+a_{n+1}) \in X^\star$, it is of the
desired type.
Otherwise all of them jump over $b$ to the right. Together with the jump $(v,w)$
this gives $F^+_b \ge n+1$, and thus $F^-_b \ge1$; that is, there is a
jump over $b$ to the left, as desired.
\end{pf}

We note that, thanks to the reflection symmetry, versions of the above
lemmas hold also for jumps to the left.

%s4.3 #&#
\subsection{General facts from Gibbs measure theory}
Here we review some techniques from the theory of Gibbs measures. Most
of these are well known; our aim is to have these presented in one
block for easier later reference.
We start with the general fact that, in the context of Gibbs measures,
energy estimates yield probabilistic bounds:

%le4.5 #&#
\begin{Lem} \label{lemprobH}
Let $\La\Subset X$, $\eta\in\Om$, $A,A' \in\F$, $c\in\R$ and $N
\ge0$.
Suppose that for every $\si\in A$ with $\si\sim_\La\eta$, we can
find a
$\si' \in A'$ with $\si' \sim_\La\eta$ such that  $H_{\La}(\si') -
H_{\La}(\si) \le c$
and every $\si'$ can be attributed to at most $N$ distinct $\si$. Then
%
%e4.14 #&#
\begin{equation}
\ga_{\La}(A|\eta) \le N \texte^{c} \ga_{\La}
\bigl(A'|\eta\bigr).
\end{equation}
In particular, for $A' := \Om$ and $N := 1$, which requires $\si
\mapsto\si'$
to be injective, we get
$\ga_{\La}(A|\eta) \le\texte^{c}$.
\end{Lem}

\begin{pf}
Using the two assumptions we get
%
%e4.15 #&#
\begin{equation}
\quad\quad\sum_{\si\in A\colon\si\sim_\La\eta} \texte^{- H_{\La}(\si)} \le
\texte^{c} \sum_{\si\in A\colon\si\sim_\La\eta} \texte^{-
H_{\La}(\si')}
\le N\texte^{c} \sum_{\si' \in A'\colon\si' \sim_\La\eta} \texte
^{- H_{\La}(\si')}.
\end{equation}
Dividing both sides by $Z_\La(\eta)$ gives a corresponding inequality
for the specifications.
\end{pf}

An estimate of specification probabilities such as the one obtained in
Lemma~\ref{lemprobH}
implies a corresponding estimate of probabilities:

%le4.6 #&#
\begin{Lem} \label{lemFatou}
Let $A,A' \in\F$, $c > 0$ and suppose that for every $\eta\in\Om$
we have $\ga_{\La}(A|\eta) \le c \ga_{\La}(A'|\eta)$ for
all sufficiently large $\La\Subset X$ (where ``large'' is allowed to
depend on $\eta$).
Then $\mu(A) \le c \mu(A')$ for all $\mu\in\G$.
\end{Lem}

\begin{pf}
By the definition of Gibbs measure, for any increasing sequence of
finite $\La_m \uparrow X$ we have
%
%e4.16 #&#
\begin{eqnarray}
  c \mu\bigl(A'\bigr) - \mu(A) &=& \lim
_{m \to\infty} E_\mu \bigl( c \ga_{\La_m}
\bigl(A'|\cdot\bigr) - \ga _{\La
_m}(A|\cdot) \bigr)\nonumber
\\[-8pt]\\[-8pt]
&\ge& E_\mu \Bigl( \liminf_{m \to\infty} \bigl( c
\ga_{\La
_m}\bigl(A'|\cdot \bigr) - \ga_{\La_m}(A|
\cdot) \bigr) \Bigr) \ge0,\nonumber   %
\end{eqnarray}
where we used Fatou's lemma to get the middle inequality.
\end{pf}

Our next observation is that Gibbs measures can be obtained as weak
limits of specifications.

%le4.7 #&#
\begin{Lem} \label{lemlocal} Let $\La,\La' \Subset X$.
For every $\F_{\La'}$-measurable function $f$, the expectation $\ga
_{\La
}(f|\eta)$
is a $\F_{\La\cup\La'}$-measurable, and thus continuous, function
of~$\eta$.
\end{Lem}

\begin{pf}
Let $\eta\in\Om$. The partition function can be written as
%
%e4.17 #&#
\begin{equation}
Z_{\La}(\eta) = \sum_{\si\colon\si\sim_\La\eta}
\texte^{-
H_{\La}(\si)} = \sum_{\si'}
\texte^{-\sum_x V(\si'(x)-x)},
\end{equation}
where the exterior sum extends over all bijections
$\si'\dvtx \La\cap\eta^{-1}(\La) \to\La\cap\eta(\La)$ and
the interior
sum is over all $x \in\La\cap\eta^{-1}(\La)$. Thus $Z_{\La}(\eta)$
depends on $\eta$ only through $\eta(\La)$ and $\eta^{-1}(\La)$, and
thus is $\F_\La$-measurable.
The measurability of $\ga_{\La}(f|\eta)$ follows similarly.
\end{pf}

%le4.8 #&#
\begin{Lem} \label{lemgibbslimit}
Let  $\eta\in\Om$ and $\La_m \Subset X$ $(m \ge1)$ be such
that $\ga
_{\La_m}(\cdot|\eta)$
converges weakly, as $m \to\infty$, to some probability measure $\mu$.
Then $\mu\in\G$.
\end{Lem}

\begin{pf}
Since the specifications are consistent, for every local event $A$
and $\La\Subset X$ we have
%
%e4.18 #&#
\begin{equation}
E_\mu\ga_{\La} (A|\cdot) = \lim_{m \to\infty}E_{\ga_{\La
_m}(\cdot
|\eta)}
\ga_{\La} (A|\cdot) = \lim_{m \to\infty}\ga_{\La_m} (A|
\eta) = \mu(A)
\end{equation}
since both $1_A$ and
$\ga_{\La}(A|\cdot)$ are local and thus continuous functions.
This proves equality of $\mu$ with the measure on the extreme left
for $A$ local; an extension to general events is unique thanks to, for
example, the $\pi$-$\lambda$ theorem.
\end{pf}

%s5 #&#
\section{Proofs: Existence}
\label{sec5}
In this section we address a.s. finiteness of the flow in all Gibbs
measures, nonexistence of Gibbs measures in the absence of
interactions and tightness of the family of specifications leading to
the proof of existence of Gibbs measures. In particular, we provide
formal proofs of Lemmas \ref{lemnonexistence}, \ref{lemexistence} and~\ref{lemfiniteflow}. We assume throughout that $X$
and $V$ satisfy (X1), (V1) and (V2); all other assumptions will be
mentioned explicitly whenever needed.

%s5.1 #&#
\subsection{Finiteness of the flow}
We begin by proving Lemma~\ref{lemfiniteflow}. First we dismiss the
issue of measurability:

\begin{pf*}{Proof of Lemma~\ref{lemfiniteflow}, tail measurability}
In order to show that $F$ is $\T$-measurable, let $a \in X^\star$ and
$\La_a := \{x \in X\dvtx  x>a\}$.
Based on its definition, $F_a^+$ is measurable with respect to $\si(\{
P^-_x\dvtx  x \in\La_a\})$ and
$F_a^-$ is measurable with respect to $\si(\{P^+_x\dvtx  x \in\La_a\}
)$. Since $F$ is, modulo a proviso when both $F_a^\pm=\infty$, the
difference $F_a^+-F_a^-$, it follows that also
$F$ is $\F_{\La_a}$-measurable. Since this holds for all $a$, we get
that $F$ is $\T$-measurable.
\end{pf*}

To show that $F$ is finite $\mu$-a.s.
for every Gibbs measure, we have to prove that long jumps are unlikely. This
follows from a suitable energy estimate for permutations having two
nested jumps.

%le5.1 #&#
\begin{Lem} \label{lemnestedjump}
Let $x,y\in X$ obey $x \le y$ and let $\mu\in\G$.
Then there is an $l_0=l_0(y-x) \ge0$ such that
%
%e5.1 #&#
\begin{equation}
\mu(x\to y, v \to w) \le\frac{1} {|w-v|^5}
\end{equation}
holds for all $v,w\in X$ that obey $v \le x-1$, $w \ge y+1$ and $w-v
\ge l_0$.
\end{Lem}

\begin{pf}
For $x,y,v,w$ as above, set $l_0 := c_\psi(y-x,3)$ where $c_\psi$ is as
in \eqref{cpsi}.
Let $\eta\in\Om_X$ and pick $\La\Subset X$ such that $x,y,v,w \in
\La$.
Let $\si\sim_\La\eta$ be such that $\si(x) = y$ and $\si(v) = w$.
By our choice of $\La$ we have $\si_{vx} \sim_\La\si$ and
the energy estimate \eqref{energyestimate} implies
%
%e5.2 #&#
\begin{eqnarray}
  &&H_{\La}(\si)- H_{\La}(
\si_{vx})\nonumber
\\
&&\quad\quad\ge2 \min\{x-v,w-y\} \psi_{y-x}(w-v)\log(w-v)
\\
&&\quad\quad\ge5 \log(w-v)  \nonumber %
\end{eqnarray}
thanks to our choice of $l_0$. We also note that $\si_{vx}$ uniquely
determines $\si$
for given $x,y,v,w$. By Lemma~\ref{lemprobH} we thus get
%
%e5.3 #&#
\begin{equation}
\ga_\La(x \to y, v \to w|\eta) \le\frac{1}{|w-v|^5}.
\end{equation}
The conclusion follows by integrating with respect to  $\mu$.
\end{pf}

%le5.2 #&#
\begin{Lem} \label{lemfiniteF}
Suppose that $a \in X^\star$ satisfies $\cg(a)<\infty$. Then
every $\mu
\in\G$ obeys $\mu(F_a^+ = \infty) = 0$.
\end{Lem}

\begin{pf}
Any $\si\in\Om_X$ such that $F_a^+(\si) = \infty$ contains a jump $(x,y)$
with $x < a < y$. Fix $l_0 \ge0$ as in the previous lemma.
In addition, for any $l>0$, $\si$ also contains a
jump $(v,w)$ such that $v \le x-1$, $w \ge y+1$ and $w-v \ge l,l_0$.
(This is because only finitely many jumps can start or end in a given
finite set.)
A union bound and Lemma~\ref{lemnestedjump} then give
%
%e5.4 #&#
\begin{eqnarray}
 \mu\bigl(F_a^+ = \infty\bigr) &\le&
\mathop{\mathop{\sum_{x,y,v,w\in X}}_{v+1
\le x < a < y \le w-1\colon}}_{
w-v \ge l,l_0}
\mu(x \to y,v \to w)\nonumber\\[-8pt]\\[-8pt]
&\le&\mathop{\mathop{\sum_{x,y,v,w\in X}}_{v <x < a < y < w\colon}}_{ w-v \ge l}
\frac{1}{|w-v|^5}. \nonumber  %
\end{eqnarray}
We first sum over $x$ and $y$ using
$\# \{x \in X\dvtx  v < x < a\} \le\cg(a) (a-v) \le\cg(a)(w-v)$
and similarly for $\{y \in X\dvtx  a < y < w\}$. In order to sum
over $v$ and $w$ we use
\eqref{sumscga2} with the result
%
%e5.5 #&#
\begin{eqnarray}
  \mu\bigl(F_a^+ = \infty\bigr) &\le&\sum
_{v <a < w\colon w-v \ge l} \frac{\cg
(a)^2}{|w-v|^3} \le\sum
_{v <a < w} \frac{\cg(a)^2}{\max\{|w-v|,l\}^3}\nonumber
\\[-8pt]\\[-8pt]
&\le&\cg(a)^4 \int_0^\infty\textd s\,
\frac{s}{\max\{s,l\}^3} = \frac
{3\cg(a)^4}{2l}.\nonumber   %
\end{eqnarray}
As $l$ can be arbitrarily large, we are done.
\end{pf}

\begin{pf*}{Proof of Lemma~\ref{lemfiniteflow}, a.s. finiteness}
Let $\mu\in\G$. Since (V1), (V2) and (X1), (X2) are assumed, all
previous derivations are at our disposal. The previous lemma gave $\mu
(F_a^+ = \infty) = 0$ and, thanks to obvious symmetry considerations
and the fact that $\cg(a)$ is defined symmetrically, we similarly have
$\mu(F_a^- = \infty) = 0$. Therefore, the flow is finite $\mu$-a.s.,
as desired.
\end{pf*}

%s5.2 #&#
\subsection{No Gibbs measures without interaction}
Here we address the fact that, in the absence of ``interactions'' the
set of Gibbs measures is empty.

\begin{pf*}{Proof of Lemma~\ref{lemnonexistence}}
Assume $V := 0$ and suppose that there is some $\mu\in\G$.
Fix $x,y \in X$, let $N \in\N$ and pick $\eta\in\Om$.
Let $\La_N \Subset X$ be such that $x,y \in\La_N$ and $\# \La_N \ge
N$. For any finite $\La\Subset X$ that obeys $\La_N \cup\eta(\La_N)
\cup\eta^{-1}(\La_N) \subset\La$ we observe that, since $\ga_\La
(\cdot
|\eta)$ is the uniform distribution on all $\si\sim_\La\eta$
and $\si
(x)$ can be any point of $\La_N$ with equal probability,
%
%e5.6 #&#
\begin{equation}
\ga_\La(x \to y|\eta) =\frac{1}{|\Lambda_N|}\le\frac{1} {N}.
\end{equation}
By Lemma~\ref{lemFatou} this implies $\mu(x \to y) \le\frac{1} {N}$
and, since $N$ was arbitrary, $\mu(x \to y) = 0$ for all $y \in X$, a
contradiction. Hence $\G=\varnothing$ after all.
\end{pf*}

%s5.3 #&#
\subsection{Existence of Gibbs measures}
Now assuming (X2) in addition to conditions (X1), (V1) and (V2),
we proceed to establish Lemma~\ref{lemexistence} dealing with existence
of Gibbs measures. Not unexpectedly, this will be done by proving
tightness for sequences of specifications over increasing volumes. We
begin by deriving an estimate on the probability of a long jump.

%le5.3 #&#
\begin{Lem} \label{lemjump}
Let $\eta\in\Om_X$ be such that $|F(\eta)| =: n \in\N$, and fix $x
\in X$.
There are $\La_0 \Subset X$
and $l_0 \in\R$ (both depending on $n,x$ only) such that for
all $\Lambda$ with
$\La_0 \subset\La\Subset X$ and all $v,w\in X$ with $v \le x \le w$ and
$w-v \ge l_0$ we have
%
%e5.7 #&#
\begin{equation}
\ga_\La(v \leftrightarrow w|\eta) \le\frac{2}{|w-v|^3}.
\end{equation}
\end{Lem}

\begin{pf}
Let $\La_0 \Subset X$ be so large that it contains all jumps of $\eta$
to, from and over any point of the set $\{x_{i}\dvtx -n \le i \le n\}$.
This is possible because $\eta$ has only finitely many jumps
over these points, in light of finiteness of $F(\eta)$.
Define, assuming (V1) and (V2),
%
%e5.8 #&#
\begin{equation}
l_0 := \max \bigl\{x_{n}-x_{-n}+1,c_\psi
\bigl(x_{n}-x_{-n},2 \cS (x,n)\bigr) \bigr\},
\end{equation}
and suppose that $v,w\in X$ are such that
$v \le x \le w$ and $w-v \ge l_0$. Note that if $\Lambda$ obeys $\La_0
\subset\La\Subset X$ and if $\si\sim_\La\eta$ is such that $\si(v)
= w$, the definition of $l_0$ implies $w > x_{n}$ or $v< x_{-n}$. We
will only consider
the first case (the other one can be done similarly).

By Lemma~\ref{lemfindjump} we can choose $y,z \in X$ such that
$\si(y) = z$ and $y \ge x_{1}$ and $z \le x_{n}$ hold.
Then $y,z,v,w\in
\Lambda_0$ and so $y,z,v,w\in\Lambda$. It follows that $\si_{vy}
\sim
_\La\si\sim_\La\eta$, and we can use
the energy estimate \eqref{energyestimate} from Lemma~\ref{lemenergyestimate}
to obtain
%
%e5.9 #&#
\begin{eqnarray}
 &&H_{\La}(\si)-  H_{\La}(
\si_{vy})\nonumber
\\
&&\quad\quad\ge2 \min\{x_{1}-v,w-x_{n}\} \psi_{x_{n}-x_{1}}(w-v)
\log(w-v)
\\
&&\quad\quad\ge\frac{2}{\cS(x,n)} \psi_{x_{n}-x_{-n}}(w-v) \log(w-v) \ge3 \log(w-v).
\nonumber  %
\end{eqnarray}
We also note that $\si_{vy}$ uniquely determines $\si$ for
given $v,w$ [since
$\si_{vy}^{-1}(w) = y$]. By Lemma~\ref{lemprobH}
we thus get $\ga_\La(v \ra w|\eta) \le\frac{1} {|w-v|^3}$, as desired,
and by symmetry
we get the same estimate for jumps $w \to v$ to the left.
\end{pf}

%le5.4 #&#
\begin{Lem} \label{lemtight}
Suppose $\eta\in\Om_X$ obeys $|F(\eta)| <\infty$, and let $\La_m
\Subset X$
be sets such that $\La_m \uparrow X$ for $m \to\infty$.
Then $\{\ga_{\La_m}(\cdot|\eta)\}_{m\ge1}$ is tight.
\end{Lem}

\begin{pf}
Let $\ep>0$ and choose $a(x) > 0$ such that $\sum_{x \in X} a(x) \le1$.
For $x \in X$ let
%
%e5.10 #&#
\begin{equation}
K_f(x) := \bigl\{\si\in\Om_X\dvtx \bigl|x-\si(x)\bigr|,\bigl|x-
\si^{-1}(x)\bigr| \le f(x) \bigr\},
\end{equation}
where $f(x) \ge0$ is chosen so that the following holds true: First we need
$f(x) \ge l_0$, where $l_0=l_0(n,x)$ with $n:=|F(\eta)|$ are as in
previous lemma. Next we will require that
%
%e5.11 #&#
\begin{equation}
f(x) \ge\sqrt{\frac{6\cg(x)}{\ep a(x)}}
\end{equation}
and also
%
%e5.12 #&#
\begin{equation}
\label{tight} \ga_{\La_m}\bigl(K_f(x)^\cc|\eta
\bigr) \le\ep a(x)
\end{equation}
for every $m$ such that $\La_m \not\supset\Lambda_0$,
where $\Lambda
_0=\La_0(x,\eta)$ is from the previous lemma.
The latter is possible since there are only finitely many such $\La_m$
and $K_f(x) \uparrow\Om_X$ for $f(x) \to\infty$.
For $m$ such that $\La_m \supset\La_0(x,\eta)$ and $w \ge x$ such
that $w-x \ge f(x)$
we can use Lemma~\ref{lemjump} to estimate $\ga_{\La_m}(x
\leftrightarrow w|\eta)
\le\frac{2} {|w-v|^3}$,
and using \eqref{sumscga} gives
%
%e5.13 #&#
\begin{eqnarray}
  \sum_{w > x+f(x)}
\frac{2}{(w-x)^3} &\le&\sum_{w > x} \frac{2}{\max\{w-x,f(x)\}^3}\nonumber
\\[-8pt]\\[-8pt]
& \le&\int_0^\infty \textd s \frac{2\cg(x)} {\max\{s,f(x)\}^3} =
\frac{3 \cg(x)}{f(x)^2}.\nonumber   %
\end{eqnarray}
By symmetry we have the same estimate for jumps from and to $x$ from
the left.
Due to the choice of $f(x)$ and the definition of $K_f(x)$, the bound
\eqref{tight} thus holds for all $m\ge1$.

Now set $K_f := \bigcap_{x \in X} K_f(x)$ and note that
%
%e5.14 #&#
\begin{equation}
\ga_{\La_m}\bigl(K_f^\cc|\eta\bigr) \le\sum
_x \ga_{\La_m}\bigl(K_f(x)^\cc
|\eta\bigr) \le \sum_x \ep a(x) \le\varepsilon
\end{equation}
for all $m\ge1$. Since $K_f$ is compact, which can be seen by a diagonal
argument using the completeness of $\Om_X$, we obtain
tightness for the sequence $\{\ga_{\La_m}(\cdot|\eta)\}_{m\ge1}$.
\end{pf}

Our last item of concern is whether or not subsequential limits of
specifications are Gibbs measures with the correct value of the flow.

%le5.5 #&#
\begin{Lem} \label{lemlimitflow}
Let $\eta\in\Om_X$ such that $n:=F(\eta)\in\Z$. Let $\La_m
\Subset
X$ be sets such that $\La_m \uparrow X$ and such that $\ga_{\La
_m}(\cdot
|\eta)$ converges, as  $m \to\infty$, weakly
to a probability measure $\mu$.
Then $\mu\in\G_n$.
\end{Lem}

\begin{pf} Lemma~\ref{lemgibbslimit} ensures $\mu\in\G$
and so we only have to show that $\mu(F_a= n) = 1$ for some $a \in
X^\star$.
For this it is convenient to define a localized
version of the flow: For a given $l \ge0$ let $F_a^{l,\pm}(\sigma)$ be
the analogues of the quantities $F_a^\pm(\sigma)$, respectively, but
counting only jumps over $a$ of length at most $l$. (The specific
choice of ``length'' is not important here, one can, e.g., use
Euclidean distance.) Define $F_a^l(\sigma):=F_a^{l,+}(\sigma
)-F_a^{l,-}(\sigma)$, where no provisos are necessary because all
quantities are finite. Although $F_a^l(\sigma)$ may,
unlike $F_a(\sigma
)$, depend on $a$, we have
%
%e5.15 #&#
\begin{equation}
\label{E:5.10} \bigl|F_a(\sigma) \bigr|<\infty \quad\Rightarrow\quad
F_a^l(\sigma ) \mathop{\longrightarrow}\limits
_{l\to\infty}
F_a(\sigma).
\end{equation}
This is because the finiteness of $F_a(\sigma)$ implies that all jumps
over $a$ have a bounded length (depending only on $\sigma$).

Returning to the main line of the proof, we note that
%
%e5.16 #&#
\begin{eqnarray}
%
  %\label{}
\label{faequalsn} \bigl|\mu(F_a=n)
- 1 \bigr| &=& \bigl|\mu(F_a=n) - \ga_{\La
_m}(F_a=n|\eta )
\bigr|\nonumber
\\
&\le&\mu\bigl(F_a \neq F_a^l\bigr) +
\ga_{\La_m}\bigl(F_a^l\neq F_a|\eta
\bigr)
\\
&&{} + \bigl|\mu\bigl(F_a^l=n\bigr) - \ga_{\La_m}
\bigl(F_a^l=n|\eta\bigr) \bigr| . \nonumber  %
\end{eqnarray}
If $F_a(\si) \neq F_a^l(\si)$, then $\si$ contains a jump to, from
or over
$x := a_{1}$ of length $\ge l$, so if $m$ and $l$ are sufficiently large
(depending on $\eta,a$ only), Lemma~\ref{lemjump} and \eqref
{sumscga} give
%
%e5.17 #&#
\begin{eqnarray}
\label{E:5.12} %
  \ga_{\La_m}
\bigl(F_a \neq F_a^l|\eta\bigr) &\le&\sum
_{v < a < w\colon w-v \ge l} \ga_{\La_m}(v \leftrightarrow w|\eta)\nonumber
\\
&\le&\sum_{v < a < w\colon w-v \ge l} \frac{2} {(w-v)^3}\nonumber
\\[-8pt]\\[-8pt]
&\le&\sum_{v < a < w} \frac{2} {\max\{l,w-v\}^3}\nonumber
\\
&\le&\int_0^\infty \textd s \,\frac{2\cg(a)^2s} {\max\{l,s\}^3} =
\frac{3 \cg(a)^2}l. \nonumber  %
\end{eqnarray}
Since $\{F_a^l = n\}$ is a local event and thus its indicator is a
continuous function,
we definitely have $\ga_{\La_m}(F_a^l=n|\eta) \to\mu(F_a^l=n)$ as  $m
\to\infty$,
so letting $m \to\infty$ in \eqref{faequalsn} yields
%
%e5.18 #&#
\begin{equation}
\bigl|\mu(F_a=n) - 1 \bigr| \le\mu\bigl(F_a \neq
F_a^l\bigr) + \frac{3 \cg(a)^2}l,
\end{equation}
once $l$ is sufficiently large. Since $\mu\in\G$, Lemma~\ref
{lemfiniteflow} implies $|F_a|<\infty$ $\mu$-a.s. and so, by \eqref
{E:5.10}, also $F_a^l\to F_a$ $\mu$-a.s. as $l\to\infty$. In
particular, $\mu(F_a \neq F_a^l)\to0$ in this limit as well and so the
claim follows by taking $l\to\infty$.
\end{pf}

\begin{pf*}{Proof of Lemma~\ref{lemexistence}}
This is now a direct consequence of the three preceding lemmas.
\end{pf*}

%s6 #&#
\section{Proofs: Infinite cycles and uniqueness}
\label{sec6}
In this final section we develop the desired level of control over the
number of infinite cycles in permuatations sampled from a Gibbs
measure. The key notion is that of a cut, introduced in Definition~\ref
{def1}. Cuts will allow us to give full classification of all Gibbs
measures leading to the proof of Theorem~\ref{thm}. Naturally, we will
also provide formal proofs of Lemmas \ref{leminfcycles} and \ref
{lemuniqueness}. Throughout this section we assume the validity of
(V1), (V2), (X1) and (X$n$), where $n$ will be clear from context.

%s6.1 #&#
\subsection{Pre-cuts}
\label{secgood}
As already mentioned, cuts permit us to control the exact number of
infinite cycles crossing over a given point of $X^\star$. In order to
exercise this control throughout $X^\star$, we need to identify a
bi-infinite sequence of cuts in a.e. permutation. This will be
achieved by relaxing to the notion of $k$-pre-cut, defined as follows:

%de6.1 #&#
\begin{definition}
For $k \in\N$, $\si\in\Om_X$ and $n:=F(\sigma)\in\Z$, the
point $a\in
X^\star$ is called a $k$-pre-cut for $\si$ if all jumps
of $\si$ over $a$, from $a_{-n},\ldots,a_{-1}$ or to $a_{1},\ldots,a_{n}$
are completely contained in $\{a_{-k},\ldots,a_{k}\}$.
\end{definition}

The role of $k$-pre-cuts is that of candidates for $k$-cuts, since
any $k$-pre-cut $a$ can
be made into a cut by modifying the permutation locally near $a$.

In this section we prepare the proof of Lemma~\ref{leminfcycles} by
showing that, in a.e. permutation,
there are sufficiently many $k$-pre-cuts. Let $n \ge0$ be fixed
throughout this
section, and let $c_n$ be a sequence sufficiently large that $c_n \ge
n+1$ and such that
%
%e6.1 #&#
\begin{equation}
X^\star_n := \bigl\{a \in X^\star\dvtx \cS(a,n)
\le c_n, \cg(a) \le c_n \bigr\},
\end{equation}
is bi-infinite. The latter is possible in light of condition (X$n$).
Let $k$ be the smallest natural number such that
%
%e6.2 #&#
\begin{equation}
\quad\quad k \ge\max \bigl\{n+1, c_n c_\psi(nc_n,2c_n),
c_nc_\psi(0,2c_n)+nc_n^2,
24 c_n^3, 33n c_n^2+1 \bigr\}.
\end{equation}
We first consider possible configurations in which a point $a \in
X^\star$ is not a
$k$-pre-cut. For this we call a jump $v \to w$ (with $v,w \in X$, of
course) $a$-relevant,
if it is a jump over $a$ that is going either from $a_{-n},\ldots
,a_{-1}$ or to $a_{1},\ldots,a_{n}$.

%le6.2 #&#
\begin{Lem} \label{lembad}
Let $\mu\in\G_n$, $a \in X^\star_{n}$ and $k$ as above. For $\De
\subset X$
let $C(\De,a)$ be the event that all $a$-relevant jumps are contained
in $\De$
and let $B \in\F_{\De^\cc}$.
Then for all $a$-relevant jumps $v \to w$ ($v,w \in X$) that are not contained
in $\{a_{-k},\ldots,a_{k}\}$ we have
%
%e6.3 #&#
\begin{equation}
\mu \bigl(C(\De,a) \cap B \cap\{v \to w\} \bigr) \le\frac{1} {|w-v|^3} \mu(B).
\end{equation}
\end{Lem}

The proof is easier to present when we deal separately with the case
when the said jump $v\to w$ is to the right and to the left.

\begin{pf*}{Proof of Lemma~\ref{lembad}, the case $w>v$}
Here we suppose that $w>v$ which means that $v \to w$ is a jump to the
right. In this case, either $v < a < a_{k}< w$ or $v<a_{-k}< a < w$.
From these two possibilities we will address only the former since the
latter is quite analogous.

Let $\eta\in\Om_X$ be such that $F(\eta) = n$.
Let $\La_0 \Subset X$ be so large that it contains all jumps of $\eta$
over $a$ or to $a_{1},\ldots,a_{n}$
[which is possible since $F(\eta)$ is finite], and let $\La\Subset X$
be such that $\La_0 \subset\La$.
Suppose $\si\sim_\La\eta$ is such that $\si(v) = w$ and $\si\in
C(\De
,a) \cap B$.
By Lemma~\ref{lemfindjump} we can find a jump $y\rightarrow z$ such that
$y \ge a_{1}$, $z \le a_{n}$. Thanks to the containment $\Lambda
_0\subset\La$ we have $v,w,y,z \in\La$,
so $\si_{vy} \sim_\La\si\sim_\La\eta$
and the energy estimate from Lemma~\ref{lemenergyestimate}
(with $y':=a_1$ and $z':=a_n$) gives
%
%e6.4 #&#
\begin{eqnarray}
 \quad H_{\La}(\si)- H_{\La}(
\si_{vy}) &\ge&2 \min\{a_{1}-v,w-a_{n}\}
\psi_{a_{n}-a_{1}}(w-v) \log(w-v)\nonumber
\\[-8pt]\\[-8pt]
&\ge&\frac{2}{c_n} \psi_{nc_n}(w-v) \log(w-v) \ge3 \log(w-v).
 \nonumber %
\end{eqnarray}
Here in the second bound we used that $a \in X^\star_n$ in order
to estimate minimal and maximal distances against $c_n$ and then
applied the natural monotonicity of $d\mapsto\psi_d(x)$. The last
inequality holds because $w-v \ge a_{k}-a \ge\frac{k} {c_n} \ge c_\psi
(nc_n,2c_n)$.

Note that for $\si\in C(\De,a)$ we have $v,w,y,z \in\De$, so $\si
\in
B$ implies \mbox{$\si_{vy} \in B$}. The swapped permutation $\si_{vy}$
uniquely determines $\si$ for given $v,w$ [since
$\si_{vy}^{-1}(w) = y$]. By Lemma~\ref{lemprobH}
we thus get
%
%e6.5 #&#
\begin{equation}
\ga_\La \bigl(C(\De,a) \cap B \cap\{v \to w\} |\eta \bigr) \le
\frac{1}{|w-v|^3} \ga_\La(B|\eta).
\end{equation}
As this holds for all $\Lambda$ large, Lemma~\ref{lemFatou} implies the
desired estimate.
\end{pf*}

\begin{pf*}{Proof of Lemma~\ref{lembad}, the case $v>w$}
Now suppose that $v \to w$ is a jump to the left, that is,
either $w \le a_{n}< a_{k}< v$ or $w < a_{-k} < a_{-n} \le v$.
We will henceforth assume the former as the latter can be dealt with similarly.

Let $\eta\in\Om_X$ such that $F(\eta) = n$.
Let $b:= \frac{1} 2 (a_{n}+a_{n+1}) \in X^\star$ and let $\La_0
\Subset
X$ be so large
that it contains all jumps of $\eta$ over $b$ [which is possible
since $F(\eta)$ is finite], and let again $\Lambda\Subset X$
obey $\La
_0 \subset\La$.
Let $\si\sim_\La\eta$ such that $\si(v) = w$ and $\si\in C(\De,a)
\cap B$.
Since $F_b(\eta) = n \ge0$ and
$(v,w)$ is a jump over $b$ to the left, there are at least $n+1$ jumps over $b$
to the right, so one of them, say $(y,z)$ satisfies $y \le a$ and $z
\ge a_{n+1}$.
By our choice of $\La$ we again have $v,w,y,z \in\La$, so $\si_{vy}
\sim_\La\si\sim_\La\eta$
and the energy estimate \eqref{energyestimate} with reversed
directions gives
%
%e6.6 #&#
\begin{eqnarray}
&&  H_{\La}(\si)- H_{\La}(
\si_{vy})\nonumber
\\
 &&\quad\quad\ge 2 \min\{v-a,a_{n+1}-w\} \psi_{a-a_{n+1}}(v-w) \log(v-w)
\\
 &&\quad\quad\ge3\log|w-v|.\nonumber  %
\end{eqnarray}
Here we have estimated the minimal distance by $\frac{1} {c_n}$ and
used
%
%e6.7 #&#
\begin{equation}
v-w \ge a_{k}-a_{n} \ge(a_{k}-a) -
(a_{n}-a) \ge\frac{k}{c_n} - nc_n \ge
c_\psi(0,2c_n),
\end{equation}
which implies $\psi_{a-a_{n+1}}(v-w) \ge\psi_0(v-w) \ge2 c_n$.
All estimates involving $c_n$ use that $a \in X^\star_n$.

Thanks to the choice of $\si\in C(\De,a)$
we have $z,y,v,w \in\De$, so $\si\in B$ implies $\si_{vy} \in B$.
Moreover, $\si_{vy}$ uniquely determines $\si$ for given $v,w$ [since
$\si_{vy}^{-1}(w) = y$]. As above, a combination of Lemmas \ref
{lemprobH} and \ref{lemFatou} then implies the desired estimate.
\end{pf*}

We can now move on to the main conclusion of this subsection:

%le6.3 #&#
\begin{Lem} \label{leminfprecut}
Let $\mu\in\G_n$ and let $Y_n \subset X_n^\star$ be bi-infinite.
Then $Y_n$ contains bi-infinitely many $k$-pre-cuts $\mu$-a.s.
\end{Lem}

\begin{pf}
The proof proceeds by a sort of renewal argument: We examine a
subsequence of points from $Y_n$ in an ordered fashion and note that
the probability of not seeing a $k$-pre-cut in the first $m$ of them
decays exponentially with $m$. We will have to do this relative to any
position (marked by an integer $N$) and for any number of consecutive
points (marked by a natural $M$).

Let us now proceed with details. Fix any $N \in\Z$, let $M \in\N$ and
consider points $a_1,\ldots,a_M \in Y_n$ be such that
$a_1 \ge N$ and
%
%e6.8 #&#
\begin{equation}
(a_m)_{-n} - (a_{m-1})_{n+k} \ge l(m)
:= c_n^2 2^{m+4} \qquad\mbox{for all } 1 < m \le M.
\end{equation}
Let $A_m$ denote the event that $a_m$ is not a $k$-pre-cut, set
$B_m := \bigcap_{i \le m} A_i$ (with $B_0 := \Om_X$)
and write $C_m$ for the set of all configurations containing
an $a_m$-relevant jump
that is not contained in $\De_{m-1} := \{(a_{m-1})_{i}\dvtx  i > n\}$
(with $\De_0 := X$). Our goal is to derive an inductive bound on $\mu(B_m)$.

We will begin by deriving a bound on $\mu(C_m)$. Fix some $m$ with $1
\le m \le M$ for the time being and simplify notation by setting $a :=
a_m$ and $a' := a_{m-1}$. We note that
%
%e6.9 #&#
\begin{eqnarray}
  \mu(C_m) &\le& \sum
_{v<a'_{n}, w>a} \mu(v \to w) + \sum_{v>a_{-n}, w<a'_{n}}
\mu(v \to w)\nonumber
\\[-8pt]\\[-8pt]
&\le&\sum_{v<a<w\colon w-v>l(m)} \frac{2}{|w-v|^3} + \sum
_{1 \le i
\le n} \sum_{w <a'_{n}}
\frac{1}{|w-a_{-i}|^3}, \nonumber  %
\end{eqnarray}
where we have estimated each probability using
Lemma~\ref{lembad} with $\De:= X$ and $B := \Om_X$ noting that each
jump considered
is $a$-relevant and not contained in $\{a_{-k},\ldots,a_{k}\}$.
In the last display the first sum can be estimated using \eqref
{sumscga} and
$a \in X^\star_n$ by
%
%e6.10 #&#
\begin{equation}
\quad\quad\sum_{v<a<w\colon w-v>l(m)} \frac{2}{|w-v|^3} \le\int
_0^\infty \frac{
c_n^2 2 s}{\max\{s,l(m)\}^3}\, \textd s \le
\frac{3 c_n^2} {l(m)} \le\frac{1}{2^{m+2}}.
\end{equation}
Similarly the second sum can be estimated by
%
%e6.11 #&#
\begin{eqnarray}
  n \sum_{w <a'_{n}}
\frac{1}{(a_{-n}-w)^3} &\le&\sum_{w <a} \frac{n}{\max\{a-nc_n-w,l(m)\}^3}
\nonumber\\
&\le&\int_0^\infty\frac{nc_n}{\max\{s-nc_n,l(m)\}^3}\, \textd s
\nonumber\\[-8pt]\\[-8pt]
&=&nc_n \biggl( \frac{l(m)+nc_n}{l(m)^3} + \frac{1}{2 l(m)^2} \biggr)
\nonumber\\
&\le&\frac{2nc_n}{l(m)^2} \le\frac{1}{2^{m+2}}.\nonumber   %
\end{eqnarray}
Combining the above estimates gives $\mu(C_m) \le\frac{1} {2^{m+1}}$.\vspace*{2pt}

Moving over to an inductive bound on $\mu(B_m)$, here we
employ $B_m\subseteq C_m\cup(B_m \cap C_m^\cc)$. The probability of
the first event has just been estimated, the second we further
decompose as
%
%e6.12 #&#
\begin{equation}
B_m \cap C_m^\cc\subseteq B_{m-1}
\cap C(\De_{m-1},a_m) \cap\bigcup \{ v \to w\},
\end{equation}
where the union is taken over all $a$-relevant jumps $v \to w$ that
not contained in $\{a_{-k},\ldots,a_{k}\}$.
Since $B_{m-1} \in\F_{\De_{m-1}^\cc}$, Lemma~\ref{lembad} implies
%
%e6.13 #&#
\begin{eqnarray}
  \mu\bigl(B_m \cap C_m^\cc
\bigr) &\le& \sum\mu \bigl(B_{m-1} \cap C(\De,a_m)
\cap\{v \to w\} \bigr)\nonumber
\\[-8pt]\\[-8pt]
&\le&\mu(B_{m-1}) \sum \frac{1} {|w-v|^3}, \nonumber
\end{eqnarray}
so it remains to bound the sum on the right (which is over pairs $v\to
w$ as specified above).

Every jump $v \to w$ considered is either a jump over $a$ of length
$|w-v| \ge|a_{\pm k}-a| \ge\frac{k} {c_n}$ or a jump from $v > a_{k}$
to $w \in\{a_{1},\ldots,a_{n}\}$ of length $v-w = v-a - (w-a) \ge
\frac{k} {c_n} - nc_n$ or a jump from $v \in\{a_{-1},\ldots,a_{-n}\}$
to $w < a_{-k}$ of length $v-w \ge v-a - (w-a) \ge
\frac{k} {c_n} - nc_n$. The corresponding contributions to the above sum
can be estimated using \eqref{sumscga}. For the first case
%
%e6.14 #&#
\begin{equation}
\sum_{v < a < w\colon|v-w| \ge k/c_n} \frac{2} {|w-v|^3} \le
c_n^2 \frac{3}{k/c_n} \le\frac{1} 8,
\end{equation}
thanks to our choice of $k$. For the second case (and similarly for the
third case)
%
%e6.15 #&#
\begin{eqnarray}
  \sum_{a < w < a_{n},v > a_{k}}
\frac{1} {|w-v|^3} &\le& n \sum_{v > a_{k}}
\frac{1} {|v-a-nc_n|^3}\nonumber
\\
&\le& n c_n \int_{0}^\infty
\frac{1} {(\max\{s, k/c_n\} - nc_n)^3
}\, \textd s\nonumber
\\[-8pt]\\[-8pt]
&=& nc_n \biggl(\frac{k/c_n}{(k/c_n-nc_n)^3} + \frac
{1}{2(k/c_n-nc_n)^2} \biggr)\nonumber
\\
&\le& nc_n \frac{2}{k/c_n-nc_n} \le\frac{1} {16},  \nonumber
\end{eqnarray}
again by our choice of $k$.

Combining the above bounds we get
$\mu(B_m \cap C_m^\cc) \le\frac{1} 4 \mu(B_{m-1})$ and thus
%
%e6.16 #&#
\begin{equation}
\mu(B_m) \le\mu(C_m) + \mu\bigl(B_m \cap
C_m^\cc\bigr) \le\frac{1} 4 \mu(B_{m-1})+
\frac{1} {2^{m+1}},
\end{equation}
for arbitrary $m \ge1$. By induction, this implies $\mu(B_M) \le
\frac{1} {2^M}$.
Since $M \ge1$ was arbitrary, the $\mu$-probability that there is no
$k$-pre-cut right of $N$ is zero. But $N$ was arbitrary (integer) too
so, as a moment's thought shows, there are infinitely many $k$-pre-cuts
in positive direction $\mu$-a.s. A completely analogous argument proves
the same for negative direction as well.
\end{pf}

%s6.2 #&#
\subsection{Cuts}
We now turn from pre-cuts to cuts. The way things are set up, every
pre-cut has a chance to be a cut with a uniformly positive probability.
This implies that a positive proportion of pre-cuts are actually cuts,
although for us it will be enough to show that the set of cuts is
bi-infinite. Throughout we keep $n \ge0$ and $k$ exactly as in the
previous section.

%le6.4 #&#
\begin{Lem}
\label{leminfcut}
Let $\mu\in\G_n$. Then $X_n^\star$ contains bi-infinitely many
cuts $\mu$-almost surely.
\end{Lem}

\begin{pf} Let $N \in\Z$ and let $Y_n(N)$ be a an infinite subset
of $\{x \in X_n^\star\dvtx  x \ge N\}$
such that between any two successive points of $Y_n(N)$ there are at
least $2k$ points of $X$.
Let $a_1 < a_2 < \cdots$ be an enumeration of $Y_n(N)$.
For $M \in\N$ and $I \subset\N$ with $|I| = M$, let $A_M(I)$ denote
the event
that $a_i, i \in I$, are the first $M$ $k$-pre-cuts in $Y_n(N)$.
For $i \in I$ let $C_i$ be the event that $a_i$ is a cut and set
%
%e6.17 #&#
\begin{equation}
B_{m} := A_{M}(I) \cap C_{1}^\cc
\cap\cdots\cap C_{m}^\cc,\quad\quad 0\le m\le M.
\end{equation}
We again aim at estimating $\mu(B_m)$ by an exponentially decaying factor.

For fixed $M$, $m$ and $I$, let $a := a_{m}$, pick
$\eta\in\Om$ with $F(\eta) = n$ and
set $\La:= [a_{-k},a_{k}]\cap X$. Let $\si\in B_{m-1}$ such
that $\si\sim_\La\eta$.
Let $\si' \sim_\La\eta$ be the permutation minimizing $H_\La$.
Since $a$ is a $k$-pre-cut with respect to  $\si$, all jumps over $a$,
to $a_{1},\ldots,a_{n}$
and from $a_{-1},\ldots,a_{-n}$ are contained in $\La$ and Lemma~\ref
{lemminenergy} implies
that $a$ is a cut in permutation $\si'$.
We note that $\si' \in B_{m-1}$ since a change on $\La$ does not affect
the $k$-pre-cut-status or cut-status of the points of $Y_n(N)$ other
than $a = a_m$.
As $H_{\La}(\si') - H_{\La}(\si) \le0$ and as there are
at most $(2k)!$ different $\si$ that give the same $\si'$
(we change only jumps contained within $\La$), Lemma~\ref{lemprobH} implies
%
%e6.18 #&#
\begin{equation}
\ga_{\La}(B_{m-1}|\eta) \le(2k)! \ga_{\La}(B_{m-1}
\cap C_{m}|\eta).
\end{equation}
Integrating with respect to  $\mu$ gives $\mu(B_{m-1}) \le(2k)! \mu
(B_{m-1} \cap C_{m})$,
which yields
%
%e6.19 #&#
\begin{equation}
\mu(B_m) = \mu\bigl(B_{m-1} \cap C_{m}^\cc
\bigr) \le \biggl(1 - \frac{1} {
(2k)!} \biggr) \mu(B_{m-1}).
\end{equation}
Inductively we thus get $\mu(B_M) \le(1 - \frac{1} {(2k)!})^M \mu(A_M(I))$.
Since $\mu$-a.s. there are infinitely many $k$-pre-cuts in $Y_n(N)$,
we can sum over all admissible $I$, and thus get that
the $\mu$-probability to have no cuts in $Y_n(N)$ is at most
$(1 - \frac{1} {(2k)!})^M$. Letting $M \to\infty$, this shows that
$Y_n(N)$ contains a cut $\mu$-a.s. Since this is true for any $N$,
the set $Y_n$ contains infinitely many cuts $\mu$-a.s. in the positive
direction;
the negative direction is then handled similarly.
\end{pf}

We are now also able to give the following:

\begin{pf*}{Proof of Lemma~\ref{leminfcycles}}
The last lemma proves part (a) for $n\ge0$; symmetry then extends this
to $n\le0$. For part (b) with $n\ge0$ it remains to note that
occurrence of a single cut restricts the number of infinite cycles
from $-\infty$ to $\infty$ to (exactly) $n$ and rules out cycles
from $\infty$ to $-\infty$ altogether. Having infinitely many such cuts
excludes cycles from $\infty$ to $\infty$ and $-\infty$ to $-\infty
$ as
well, so the statement in (b) follows. The case $n\le0$ is completely
analogous by symmetry.
\end{pf*}

%s6.3 #&#
\subsection{Uniqueness of Gibbs measures}
The sole purpose of this subsection is to give:

\begin{pf*}{Proof of Lemma~\ref{lemuniqueness}}
Let $\mu,\tilde{\mu} \in\G_n$ for $n \ge0$; the case of negative $n$
is handled by symmetry. Our principal observation is that Lemmas \ref
{leminfprecut} and \ref{leminfcut} from the preceding subsections
still hold when ``$\mu$-a.s.'' is replaced by ``$\mu\otimes\tilde
{\mu
}$-a.s.,'' and $k$ is increased somewhat.

To see this for Lemma~\ref{leminfprecut} we first note that $\mu
\otimes\tilde{\mu}$ is a Gibbs measure with respect to the product
specification
$\ga_{\La}(\cdot|\eta,\tilde{\eta}) := \ga_{\La}(\cdot|\eta)
\otimes\ga
_{\La}(\cdot|\tilde{\eta})$.
Also, $a$ is (defined to be) a pre-cut, respectively, cut in $(\si
,\tilde{\si})$ if and only if it is a pre-cut, respectively, cut
in both $\si$ and $\tilde{\si}$. Thus the only change to be made
in the proof is that in all probability estimates
for bad jumps we have to consider two cases: the jump is bad either for
$\si$ or for $\tilde{\si}$. This leads to an additional factor of 2 in
all probability estimates.
By increasing $k$ and $l(m)$ accordingly, this factor can be easily absorbed.

Concerning Lemma~\ref{leminfcut}, the only change required to the proof
is that
$(2k)!$ has to be replaced by $[(2k)!]^2$ since we now have to take
into account
all possible local rearrangements of $\si$ and $\tilde{\si}$ that make
a given
$k$-precut into a cut. This does not affect the argument and, in particular,
part (a) of Lemma~\ref{lemuniqueness} thus holds.

It remains to show that $\mu= \tilde{\mu}$. Consider a cylinder event
$A \in\F_{[-N,N]}$ for some $N \ge1$. Let $C(a,b)$ be the event that
$a$ is the last cut before $-N$ and $b$ the first cut after $N$.
[Because of (a) these cuts exist a.s.] Let $\La:= \{x \in X\dvtx  a <
x < b\}$.
We note that
%
%e6.20 #&#
\begin{equation}
\ga_{\La} \bigl( (A\times\Om) \cap C(a,b) |\eta,\tilde{\eta } \bigr) =
\ga_{\La} \bigl( (\Om\times A) \cap C(a,b) |\eta,\tilde{\eta } \bigr)
\end{equation}
for all $\eta,\tilde{\eta} \in\Om$. Indeed, if $a$ or $b$ are not cuts
with respect to
$(\eta,\tilde{\eta})$, then both sides are $0$; otherwise both boundary
conditions
can be replaced by $\tau_n$ without changing the probabilities and then
the equality follows from the fact that the product measure
$\ga_{\La}(\cdot|\tau_n,\tau_n)$ and the event $C(a,b)$ are invariant
under $(\si,\tilde{\si}) \mapsto(\tilde{\si},\si)$.
Integrating w.r.t.  $\mu\otimes\tilde{\mu}$, we now get
%
%e6.21 #&#
\begin{eqnarray}
  \mu(A) &=& \mu\otimes\tilde{\mu} (A \times\Om) = \sum
_{a,b} \mu\otimes\tilde{\mu} \bigl( (A\times\Om) \cap
C(a,b) \bigr)\nonumber
\\[-8pt]\\[-8pt]
&=& \sum_{a,b} \mu\otimes\tilde{\mu} \bigl( (\Om
\times A) \cap C(a,b) \bigr) = \mu\otimes\tilde{\mu}(\Om\times A) = \tilde{
\mu}(A).  \nonumber %
\end{eqnarray}
Since $\bigcup_{N \in\N} \F_{[-N,N]}$ forms a $\cap$-stable
generator of the sigma algebra $\F$
the above implies $\mu= \tilde{\mu}$.
\end{pf*}

%s6.4 #&#
\subsection{Classification of Gibbs measures}
It remains to formally present the proof of our main result:

\begin{pf*}{Proof of Theorem~\ref{thm}}
All conclusions except (c) follow already from the preceding lemmas.
To get also (c), let $\La_N \uparrow X$
be an increasing sequence of sets and, given $n\in\Z$, pick $\eta\in
\Om$ with $F(\eta) = n$.
Combining Lemmas \ref{lemtight}, \ref{lemlimitflow} and part (a) of the
theorem we
see that every subsequence of $\{\ga_{\La_N}(\cdot|\eta)\dvtx  N
\ge1\}
$ has a subsequence
converging to $\mu_n$.
By a standard argument, this implies $\ga_{\La_N}(\cdot|\eta) \to
\mu
_n$ weakly as $n\to\infty$. [Otherwise one could find a local event $A$
and a subsequence $N_k$ such that $\ga_{\La_{N_k}}(A|\eta)$ stays away
from ${\mu_n}(A)$ by a positive factor.]
\end{pf*}

\section*{Acknowledgments}
We appreciate the constructive remarks of the referees.

%suskaldyti doi

% imsref loaded by arune.pranskunaite, 2014-04-29 10:12:15
%

\printaddresses

\end{document}